\documentclass{amsart}
\usepackage{umlaut}
\usepackage{a4wide}
\usepackage{epsfig}
\usepackage[T1]{fontenc}    
\usepackage{ae,aecompl}     
\usepackage{amsmath}        
\usepackage{amssymb}
\usepackage{amsthm}
\usepackage{amsfonts}
\usepackage{amscd}
\usepackage{mathrsfs}       
\usepackage{bbm}            
\usepackage{bm}             
\usepackage{url}            
\usepackage{paralist}       
\usepackage{graphics}       
\usepackage{epic}          
\usepackage[abs]{overpic}   
\usepackage{psfrag}
\usepackage{xspace}         

\input xypic
\xyoption{all}
\xyoption{arc}

\newtheorem{thm}{Theorem}[section]
\newtheorem{cor}[thm]{Corollary}
\newtheorem{lem}[thm]{Lemma}
\newtheorem{prop}[thm]{Proposition}
\newtheorem{rem}[thm]{Remark}
\newtheorem{defn}[thm]{Definition}

\def\Do{{\mathsf D}}      
\def\Up{{\mathsf U}}      
\def\ACG{{\textnormal A}} 
\def\BCG{{\textnormal B}} 
\def\ADG{{\mathscr A}}    
\def\BDG{{\mathscr B}}    
\def\SOT{{\mathcal T}}    
\def\Lv{{\mathsf L}}      
\def\Rv{{\mathsf R}}      
\def\Ass{{\mathsf{Asso}}}  
\def\Perm{{\mathsf{Perm}}}

\usepackage{realisationmacros}

\begin{document}

\title{Realizations of the associahedron and cyclohedron}

\author{Christophe  Hohlweg}
\address[Christophe Hohlweg]{The Fields Institute\\
222 College Street\\
Toronto, Ontario, M5T 3J1\\ CANADA}
\email{chohlweg@fields.utoronto.ca}
\urladdr{http://www.fields.utoronto.ca/\~{}chohlweg}

\author{Carsten Lange}
\address[Carsten Lange]{Freie Universit{\"a}t Berlin\\
FB Mathematik und Informatik\\
Arnimallee 3\\
14195 Berlin,\\ GERMANY}
\email{lange@math.tu-berlin.de}

\date{September 14, 2006}

\thanks{Hohlweg is partially supported by CRC.
 This work was done during the time both authors
were at the Institut Mittag-Leffler, Djursholm, Sweden.}

\begin{abstract}
\noindent 
We describe many different realizations with integer coordinates
for the associahedron (i.e. the Stasheff polytope) and for the cyclohedron
(i.e. the Bott-Taubes polytope) and compare them to the permutahedron of 
type~$\ACG$ and~$\BCG$ respectively.

The coordinates are obtained by an algorithm which uses an oriented
Coxeter graph of type~$\ACG_{n}$ or~$\BCG_{n}$ as only input data 
and which specializes to a procedure presented by J.-L. Loday 
for a certain orientation of~$\ACG_{n}$. The described realizations have
cambrian fans of type~$\ACG$ and~$\BCG$ as normal fan. This settles a conjecture
of N.~Reading for cambrian lattices of these types. 
\end{abstract}

\maketitle

\section{Introduction}\label{sec_intro}

The associahedron $\Ass(\ACG_{n-1})$ was discovered by J.~Stasheff in 1963~\cite{stasheff} and is 
of great importance in the theory of operads. It is a simple $(n-1)$-dimensional convex polytope 
whose $1$-skeleton is isomorphic to the undirected Hasse diagram of the Tamari lattice on the 
set~$\SOT_{n+2}$ of triangulations of 
an $(n+2)$-gon (see for instance~\cite{lee}) and therefore a fundamental example of a secondary
polytope as described in~\cite{Gelfend_Kapranov_Zelevinsky}. Numerous realizations of the associahedron 
have been given, see~\cite{chapoton_fomin_zelevinsky,loday} and the references therein.

An elegant and simple realization of the associahedron by picking some of 
the inequalities for the permutahedron is due to S.~Shnider \&~S.~Sternberg~\cite{shnider_sternberg} 
(for a corrected version consider J.~Stasheff \&~S.~Shnider~\cite[Appendix~B]{stasheff2}).
The removed inequalities are related to the vertices by a well-known surjection from~$S_{n}$ 
to the set~$Y_{n}$ of planar binary trees that relates the weak order of~$\ACG_{n-1}$ with the 
Tamari lattice as described in~\cite[Sec.~9]{bjoerner_wachs}. J.-L.~Loday presented recently an 
algorithm to compute the coordinates of this realization,~\cite{loday}: label the vertices of the 
associahedron by planar binary trees with~$n+2$ leaves and apply a simple algorithm on trees to obtain
integer coordinates in~$\R^{n}$.

The associahedron fits, up to combinatorial equivalence, in a larger family of polytopes discovered
by S.~Fomin and A.~Zelevinsky~\cite{fomin_zelevinsky} (and realized as convex polytopes by
F.~Chapoton, S.~Fomin, and A.~Zelevinsky~\cite{chapoton_fomin_zelevinsky}) that is indexed by
the elements in the Cartan-Killing classification. Among these \emph{generalized associahedra}, 
the cyclohedron~$\Ass(\BCG_n)$ was first described by R.~Bott and C.~Taubes in 1994~\cite{bott_taubes} 
in connection with knot theory and rediscovered independently by R.~Simion~\cite{simion}. It is a 
simple $n$-dimensional convex polytope whose vertices are given by the set~$\SOT_{n+2}^{B}$ of
centrally symmetric triangulations of a $(2n+2)$-gon. Various realizations have also been given
in~\cite{chapoton_fomin_zelevinsky, devadoss, markl,reiner,simion}, but none of them is similar 
to Loday's realization.  

It is a natural question to ask for a construction similar to Loday's for the cyclohedron and we present such 
a construction in this article: Starting with a realization of the permutahedron~$\Perm(\BCG_{n})$, 
i.e. the convex hull of the orbit of the point~$(1,2, \ldots , 2n)$ with respect to the action of the hyperoctahedral
group, we give an explicit description of realizations of the cyclohedron by removing facets of~$\Perm(\BCG_{n})$.
Moreover we introduce an algorithm to obtain (integer) coordinates for the vertices of these
realizations.

It should also be mentioned that the associahedron and cyclohedron fit into another large family of 
polytopes, the \emph{graph associahedra} introduced by M.~Carr and S.~Devadoss~\cite{carr_devadoss}
and by M.~Davis, T.~Januszkiewicz, and R.~Scott~\cite{fundamental_groups_and_blow_ups} in the study of real
blow-ups of projective hyperplane arrangements. Generalized associahedra and graph associahedra fit into 
the class of \emph{generalized permutahedra} of A. Postnikov~\cite{postnikov} where the right-hand 
sides for the facet inequalities of the permutahedron~$\Perm(\ACG_{n-1})$ are altered. In fact, the 
associahedron and cyclohedron can be obtained from this permutahedron by changing the right-hand side 
of some facet inequalities as described for example by Postnikov~\cite{postnikov}. This description of 
the cyclohedron is obtained from ``cyclic intervals of $[n]$'', that is, the cyclohedron is seen 
as ``graph associahedron of a cycle''. 
On the contrary, the realizations given in this article see the associahedron and cyclohedron related to 
the Coxeter graph of type~$A$ and~$B$. The associahedron and cyclohedron may be obtained in many ways by 
omission of some inequalities. We explicitly describe possible choices for these facet inequalities related 
oriented Coxeter graphs of type~$A$ and~$B$ and the resulting coordinates. 
Moreover, the presented realizations of the cyclohedron are the first explicit realizations as a ``generalized 
type~$B$ permutahedron'', i.e. the cyclohedron is obtained from the permutahedron of type~$B$ by changing
the right hand side of some inequalities.

It turns out that Loday's construction is generalized in two ways
by our algorithm: Firstly, for a certain orientation of the Coxeter graph of type~$A$ our algorithm coincides
with his construction and secondly, our algorithm also works for the type B associahedron, it yields coordinates
for the cyclohedron for any oriented Coxeter graph of type~$B$.

In \S\ref{sse:IntroAsso} and \S\ref{sse:IntroCyclo} we explain our algorithm and realizations of the 
associahedron and of the cyclohedron and state the main results. These results are proved in 
\S\ref{sec:ProofAss} and \S\ref{sec:RealizCyclo} by explicitly stating a H-representation\footnote{There are two ways
to realize a polytope: The H-representation is the intersection of a finite number of closed half spaces and
the V-description is the convex hull of a finite number of points, see Ziegler~\cite{ziegler} for further deatils.} for each realization. 
Finally in \S\ref{sec:Remarks} we comment some observations concerning isometry classes and barycenters 
of these realizations. Moreover, in this section we show that the normal fans of the realizations we 
provide coincide with the cambrian fans of type~$A$ and~$B$. We settle therefore Conjecture~1.1 of 
N.~Reading~\cite{reading} in type~$A$ and~$B$.\\

\noindent {\bf N.B.} We remark that our construction yields polytopes with the combinatorial type of the 
generalized associahedra of type~$A_{n}$ and~$B_{n}$. It might be worth to mention that the combinatorics 
of all the polytopes involved is determined by their $1$-skeleton or graph, since these polytopes are simple. 
This was shown by R.~Blind and P.~Mani-Levitska,~\cite{blind_mani}, as well as by G.~Kalai,~\cite{kalai}.

\subsection{Realizations of the associahedron}\label{sse:IntroAsso}

Let $S_n$ be the symmetric group acting on the set $[n]=\{1,\dots,n\}$. As a Coxeter group of
type~$A_{n-1}$, $S_n$ is generated by the simple transpositions $\tau_i = (i,i+1)$,~$i\in [n-1]$.
The Coxeter graph $\ACG_{n-1}$ is then
\begin{figure}[h]
      \psfrag{t1}{$\tau_{1}$}
      \psfrag{t2}{$\tau_{2}$}
      \psfrag{t3}{$\tau_{3}$}
      \psfrag{tn}{$\tau_{n-1}$}
      \psfrag{dots}{$\ldots$}
      \begin{center}
            \includegraphics[width=8cm]{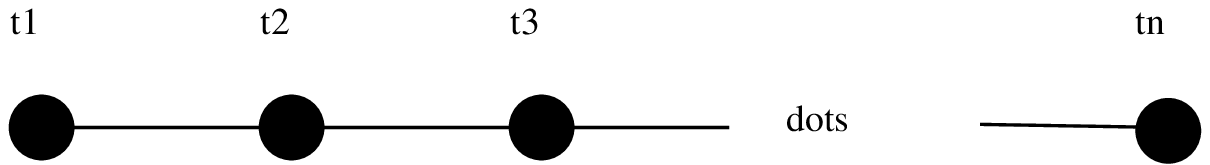}
      \end{center}
\end{figure}

Let $\ADG$ be an orientation of~$\ACG_{n-1}$. Following N.~Reading,
we distinguish between \emph{up} and \emph{down} elements 
of~$\{2, \ldots , n-1 \}$: An element $i \in \{ 2, \ldots n-1\}$ is \emph{up}
if the edge~$\{ \tau_{i-1}, \tau_{i} \}$ is directed from
$\tau_{i}$ to $\tau_{i-1}$ and \emph{down} otherwise. We extend
this definition to $[n]$ by the convention that~$1$ and~$n$ are
always down.  We remark that in Reading's work~$1$ and~$n$ can be chosen 
to be up or down. Let $\Do_\ADG$ be the set of down elements and
$\Up_\ADG$ be the set of up elements (possibly empty). 

\begin{figure}[t]
      \psfrag{0}{$0$}
      \psfrag{1}{$1$}
      \psfrag{2}{$2$}
      \psfrag{3}{$3$}
      \psfrag{4}{$4$}
      \psfrag{5}{$5$}
      \psfrag{6}{$6$}
      \psfrag{s1}{$\tau_{1}$}
      \psfrag{s2}{$\tau_{2}$}
      \psfrag{s3}{$\tau_{3}$}
      \psfrag{s4}{$\tau_{4}$}
      \begin{center}
      \begin{minipage}{0.95\linewidth}
         \begin{center}
            \includegraphics[height=6cm]{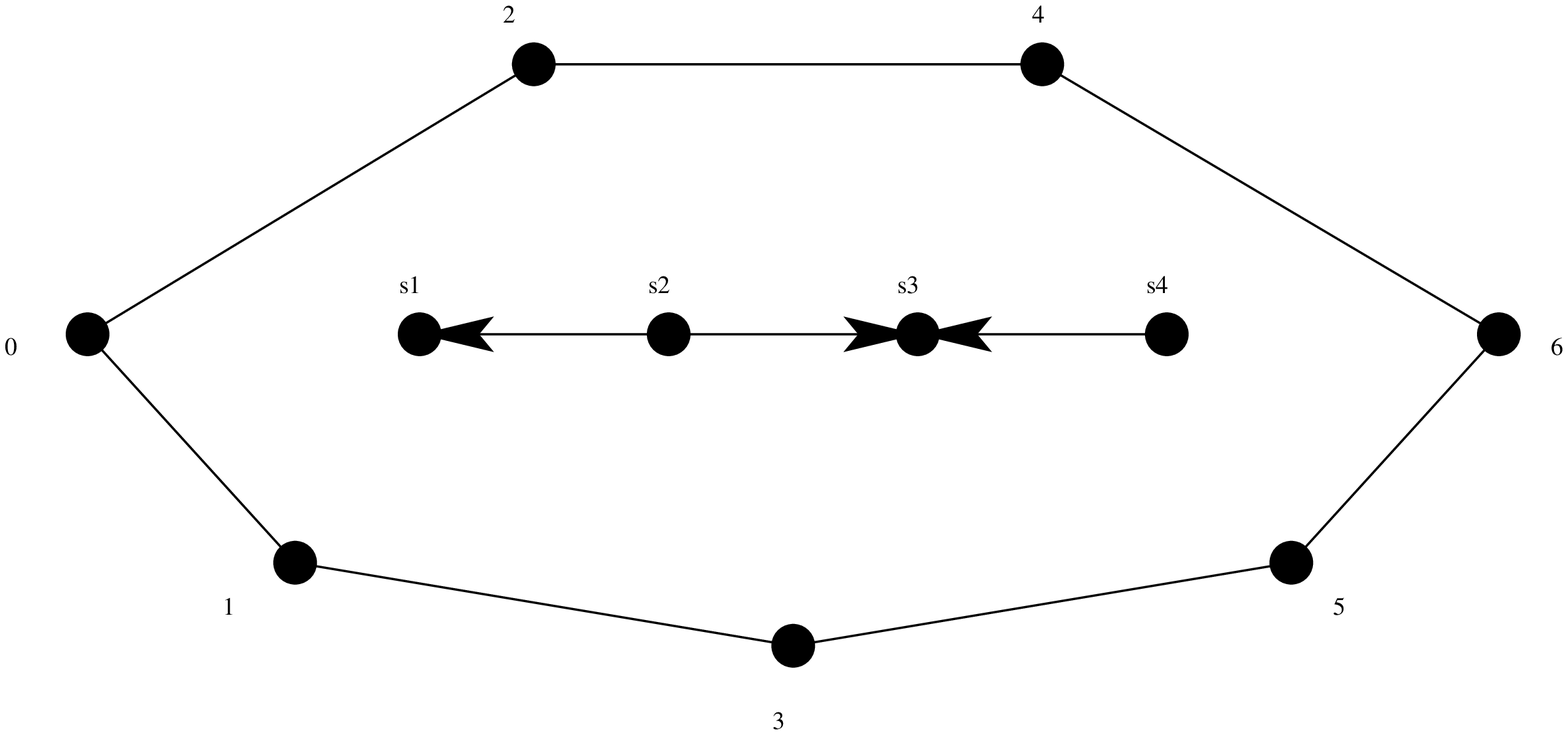}
         \end{center}
         \caption[]{A labelling of a heptagon that corresponds to
                    the orientation $\ADG$ on $\ACG_{4}$ shown inside
                    the heptagon. We have $\Do_{\ADG} = \{1, 3, 5 \}$ and
                     $\Up_{\ADG} = \{ 2,4\}$.}
         \label{fig:example_labelling}
      \end{minipage}
      \end{center}
\end{figure}

The notion of up and down induces a labelling of a {\em fixed convex}  $(n+2)$-gon $P$ as
follows: label one vertex of the $(n+2)$-gon by $0$. The adjacent
vertex in counterclockwise direction is labelled by the smallest
down element of $[n]$ not already used. Repeat this procedure as
long as there is a down element that has not been used. If there
is no such element use label $n+1$ and continue to label the next
counterclockwise element by the largest up element of $[n]$ that
has not been used so far and iterate.
An example is given in Figure~\ref{fig:example_labelling}.

Now we consider  $P$ labelled according to a fixed orientation~$\ADG$.
A triangulation~$T$ of~$P$ is a planar graph with vertex set the vertices of~$P$ and
edge set the edges of~$P$ and $n-1$ non-crossing \emph{diagonals}
different from the boundary edges. We denote by~$\SOT_{n+2}$ 
the set of all triangulations of $P$ and describe a triangulation by its non-crossing diagonals.
Our goal is now to define an injective map
\begin{align*}
   M_{\ADG}: \SOT_{n+2} &\longrightarrow  \R^n \\
                   T          &\longmapsto      (x_1,x_2,\dots,x_n)
\end{align*}
that assigns explicit coordinates to a triangulation.

Before we define~$M_{\ADG}$, we introduce a family of functions~$\mu_{i}:\{0,1,\dots,n+1\}\to [n+2]$ 
that measure distances between labels of~$P$ and that is parameterized by~$i \in [n]$. 
For~$j<i$, $\mu_i(j)$ counts the number of edges~$\{ a,b \}$ of the path (on the boundary of~$P$) connecting~$i$ and~$j$ 
that uses only labels~$\leq i$. For~$j\geq i$, $\mu_i(j)$ counts the number of edges~$\{ a,b \}$ of
the path (on the boundary of~$P$) connecting~$i$ and~$j$ that uses only labels~$\geq i$. For instance, we have $\mu_4(5)=2$ 
and $\mu_5(4)=5$ in Figure~\ref{fig:example_labelling}.
\begin{figure}[b]
      \begin{center}
      \begin{minipage}{0.95\linewidth}
         \begin{center}
         \begin{overpic}
            [width=0.7\linewidth]{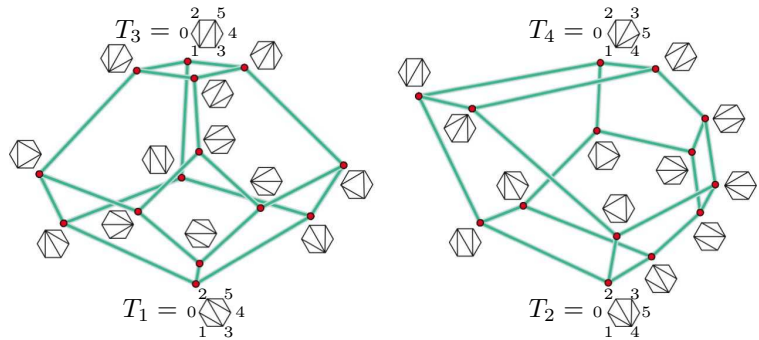}
            \put(15,1.3){$T_{1}=\ $\tiny $0$}
            \put(25,-1.3){\tiny $1$}
            \put(25,3.7){\tiny $2$}
            \put(28.5,-1.3){\tiny $3$}
            \put(30,1.3){\tiny $4$}
            \put(28.5,3.7){\tiny $5$}
            \put(14,38.4){$T_{3}=\ $\tiny $0$}
            \put(24,36){\tiny $1$}
            \put(24,41.1){\tiny $2$}
            \put(27.5,36){\tiny $3$}
            \put(29,38.4){\tiny $4$}
            \put(27.5,41.1){\tiny $5$}
            \put(69,1.3){$T_{2}=\ $\tiny $0$}
            \put(79,-1.3){\tiny $1$}
            \put(79,3.7){\tiny $2$}
            \put(82.5,3.7){\tiny $3$}
            \put(82.5,-1.3){\tiny $4$}
            \put(84,1.3){\tiny $5$}
            \put(69,38.4){$T_{4}=\ $\tiny $0$}
            \put(79,36){\tiny $1$}
            \put(79,41.1){\tiny $2$}
            \put(82.5,41.1){\tiny $3$}
            \put(82.5,36){\tiny $4$}
            \put(84,38.41){\tiny $5$}
         \end{overpic}
         \end{center}
         \caption[]{The vertices of the two associahedra shown have coordinates that are computed from 
                    triangulations of labelled hexagons.}
         \label{fig:a3_associahedra}
      \end{minipage}
      \end{center}
\end{figure}
For a triangulation~$T\in \SOT_{n+2}$ of~$P$ and $i \in [n]$, we denote by~$\Lv^T_i$ the set 
$\set{a}{0 \leq a < i \text{ and } \{a,i\} \in T}$
and by~$\Rv^T_i$ the set $\set{b}{i < b \leq n+1 \text{ and } \{b,i\} \in T}$. Set
\[
 p_{\ell}^{T}(i):=\max_{a \in \Lv_i^T}\ \{\mu_i(a)\}
 \qquad\textrm{and}\qquad
 p_{r}^{T}(i):=\max_{b \in \Rv_i^T}\  \{\mu_i(b)\}.
\]
The {\em weight}~$\omega_{i}$ of~$i$ in~$T$ is the integer~$p_{\ell}^{T}(i) p_{r}^{T}(i)$.
We now define the coordinates~$x_{i}$ of~$M_{\ADG}(T)$:
\[
  x_i := \begin{cases}
            \omega_i     & \textrm{if } i\in\Do_\ADG\\
            n+1-\omega_i & \textrm{if } i\in\Up_\ADG.
         \end{cases}
\]
In the setting of Figure~\ref{fig:a3_associahedra}, let~$\ADG_{1}$ denote the orientation that
yields the realization on the left and~$\ADG_{2}$ denote the orientation that yields the
realization on the right. Consider the triangulations $T_{1}=\{ \{0,3\}, \{2,3\},\{2,4\}\}$,
$T_{2}=\{ \{0,4\}, \{2,4\},\{3,4\}\}$, $T_{3}=\{ \{1,2\}, \{1,5\},\{3,5\}\}$, and
$T_{4}=\{ \{1,2\}, \{1,3\},\{1,5\}\}$ that are given by their set of diagonals
($T_1,T_3$ are triangulations of the left hexagon while $T_2,T_4$ are triangulations of the right 
hexagon). Then
\[
 M_{\ADG_{1}}(T_{1})=M_{\ADG_{2}}(T_{2})=(1,2,3,4)
 \qquad\text{and}\qquad
 M_{\ADG_{1}}(T_{3})=M_{\ADG_{2}}(T_{4})=(4,3,2,1).
\]

\begin{thm}\label{thm:AssA}
   Fix an orientation $\ADG$ on $\ACG_{n-1}$. The convex hull
   of~$\set{M_{\ADG}(T)}{T\in \SOT_{n+2}}$ is a realization of
   the associahedron~$\Ass(\ACG_{n-1})$ with integer coordinates.
\end{thm}

\noindent
This V-representation of~$\Ass(\ACG_{n-1})$ as convex hull of vertices is proved in
Section~\ref{sec:ProofAss}.

\begin{rem}\label{rem:AssA}\textnormal{
If all edges of~$\ACG_{n-1}$ are directed from left to right, then the realization 
just described coincides with the one given by Loday. In this situation, 
$\Up_{\ADG} = \varnothing$. Let $T\in \SOT_{n+2}$ and for each $i\in[n]$, let~$a$ 
and~$b$ be such that $p_\ell^T(i)=\mu_i(a)$ and $p_r^T(i)=\mu_i(b)$. Consider the 
triangle~$\{a,i,b\}$. Label this triangle by~$i$. Now, take the dual graph of~$T$: 
it is a planar binary tree with~$n+1$ leaves whose root is determined by the 
edge~$\{0,n+1\}$ of~$T$ and whose internal nodes are labelled by the label of the 
corresponding triangle. Then for each $i\in [n]$ the weight~$\omega_i$ of~$i$ is 
the product of the leaves of the left side of~$i$ and of the leaves of the right 
side. That is precisely how J.-L.~Loday computes the coordinates of the vertices 
in his realization, starting from planar binary trees.   }
\end{rem}

\noindent
We are now interested in a surjective map
\[
  \Phi_\ADG: S_n \to \SOT_{n+2}
\]
for any given orientation~$\ADG$ of~$\ACG_{n-1}$. These maps have been used earlier.
L.~Billera \&~B.~Sturmfels used them in~\cite{billera_sturmfels} where associahedra are
described as iterated fibre polytopes.  Other references are V.~Reiner~\cite{reiner}, 
A.~Bj\"orner \&~M.~Wachs~\cite[Remark~9.14]{bjoerner_wachs}, J.-L.~Loday \&~M.~Ronco,~\cite{LR}, 
A.~Tonk,~\cite{tonk}, and N.~Reading~\cite{reading}, who used this family of maps in his
study of \emph{cambrian lattices} of type~$A$. The orientation~$\ADG$ where each 
edge of~$A_{n-1}$ is oriented from left to right yields a well-studied map that turns out to
be a lattice epimorphism from the (right) weak order lattice to the Tamari lattice.
In fact, the undirected Hasse diagram of a cambrian 
lattice of type~$A$ is combinatorially equivalent to the $1$-skeleton 
of~$\Ass(\ACG_{n-1})$~\cite[Theorem~1.3]{reading}. In other words, these 
maps can be viewed as `$1$-skeleton' maps from~$\Perm(\ACG_{n-1})$ 
to~$\Ass(\ACG_{n-1})$.

\begin{figure}
      \psfrag{0}{$0$}
      \psfrag{1}{$1$}
      \psfrag{2}{$2$}
      \psfrag{3}{$3$}
      \psfrag{4}{$4$}
      \psfrag{5}{$5$}
      \psfrag{6}{$6$}
      \begin{center}
      \begin{minipage}{\linewidth}
         \begin{center}
            \includegraphics[height=5.5cm]{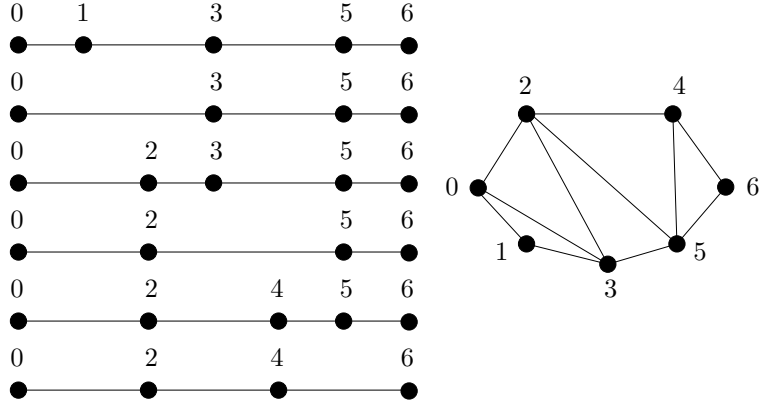}
         \end{center}
         \caption[]{We assume the same orientation $\ADG$ of $\ACG_5$ as in
                     Figure~\ref{fig:example_labelling}. The
                    permutation~$\sigma = 12345$ yields the
                    six paths shown on the left. The edges of these paths
                    form the diagonals of the triangulation~$\Phi_{\ADG}(\sigma)$
                    shown on the right.}
         \label{fig:example_triangulation}
      \end{minipage}
      \end{center}
\end{figure}

We follow the procedure given in~\cite{reading} by N.~Reading to describe 
these maps~$\Phi_{\ADG}$. We remark that N.~Reading uses the left weak order 
while we prefer the right weak order. So we invert~$\sigma \in S_{n}$ to 
translate between left and right weak order. Let $\sigma \in S_n$ and start with
the path of the labelled $(n+2)$-gon that connects~$0$ with~$n+1$
and passes through all down elements. Now read the permutation $\sigma^{-1}$
(represented as a word in $[n]$) from left to right and construct inductively
a new path from~$0$ to~$n+1$: If the next letter of $\sigma^{-1}$ is a
down element then delete this element in the path; if the next
letter is an up element then insert this element between its
largest predecessor and its smallest successor in the path. The
edges used during this process define a
triangulation~$\Phi_\ADG(\sigma)$ of the labelled $(n+2)$-gon,
see Figure~\ref{fig:example_triangulation} for an example.

\medskip
The permutahedron~$\Perm(\ACG_{n-1})$ is the classical permutahedron~$\Pi_{n-1}$
which is defined as the convex hull of the points
\[
  M(\sigma):=(\sigma(1),\sigma(2),\dots,\sigma(n))\in \R^n,\qquad \forall \sigma\in S_n .
\]
The idea of~S.~Shnider \&~S.~Sternberg to obtain the associahedron from the permutahedron by discarding inequalities
extends to all realizations of the associahedron of Theorem~\ref{thm:AssA}. The map $K_{\ADG}$ assigns
subsets of~$[n-1]$ to a diagonal and is defined in Section~\ref{sec:ProofAss}.
\begin{prop}\label{prop:AssAfromPermA}
   Fix an orientation~$\ADG$. The associahedron of Theorem~\ref{thm:AssA} is given by a subset of the
   inequalities for the permutahedron~$\Perm(\ACG_{n-1})$. These inequalities are determined by the
   image under~$K_{\ADG}$ of the diagonals of the $(n+2)$-gon labelled according to~$\ADG$.
\end{prop}

Moreover, the following analog to \cite[Proposition 2]{loday} shows that our realizations are closely related
to the maps~$\Phi_{\ADG}$: A triangulation~$T$ has a singleton $\{ \sigma \} = (\Phi_{\ADG})^{-1}(T)$
as preimage if and only if $M_{\ADG}(T)$ is a vertex of the permutahedron~$\Perm(\ACG_{n-1})$.
But there are more ways to characterize the common vertices of the associahedron and permutahedron which
depend on a chosen orientation~$\ADG$. 
\begin{prop}\label{prop:AssAandPermA}
   Fix an orientation $\ADG$ on $\ACG_{n-1}$ and let $T\in \SOT_{n+2}$ and $\sigma\in S_n$. 
   The following statements are equivalent:
   \begin{compactenum}[(a)]
      \item $M_{\ADG}(T) = M(\sigma)$,
      \item $\Phi_{\ADG}(\sigma)=T$ and the diagonals of~$T$ can be labelled such that
            \[
              \varnothing \subset K_{\ADG}(D_1) \subset \ldots \subset K_{\ADG}(D_{n-1}) \subset [n]
            \]
            is a sequence of strictly increasing nested sets.
      \item $(\Phi_{\ADG})^{-1}(T)=\{\sigma\}$,
      \item $\Phi_{\ADG}(\sigma) = T$ and for each~$i\in[n]$ we have~$p_{\ell}^{T}(i)=1$ or~$p_{r}^{T}(i)=1$.
   \end{compactenum}
\end{prop}

\noindent
The proof of both propositions is postponed to Section~\ref{sec:ProofAss}.

\subsection{Realizations of the cyclohedron}\label{sse:IntroCyclo}

An orientation~$\ADG$ of~$\ACG_{2n-1}$ is {\em symmetric} if the edges $\{\tau_i,\tau_{i+1}\}$ 
and $\{\tau_{2n-i-1},\tau_{2n-i}\}$ are oriented in \emph{opposite directions} for all~$i\in [2n-2]$. 
There is a bijection between the symmetric orientations of~$\ACG_{2n-1}$ and the orientations~$\BDG$ 
of the Coxeter graph~$\BCG_{n}$ that we describe below.  A triangulation $T\in \SOT_{2n+2}$ is 
{\em centrally symmetric} if~$T$, viewed as a triangulation  of a {\em regular} $(2n+2)$-gon, 
is centrally symmetric. Let $\SOT_{2n+2}^B$ be the set of the centrally symmetric triangulations 
of the labelled $(2n+2)$-gon.

\begin{thm}\label{thm:AssB}
   Let~$\ADG$ be a symmetric orientation of $\ACG_{2n-1}$. The convex hull
   of $\set{M_{\ADG}(T)}{T\in \SOT^B_{2n+2}}$ is a realization of the
   cyclohedron~$\Ass(\BCG_{n})$ with integer coordinates.
\end{thm}

\begin{figure}[t]
      \begin{center}
      \begin{minipage}{0.95\linewidth}
         \begin{center}
            \includegraphics[width=0.8\linewidth]{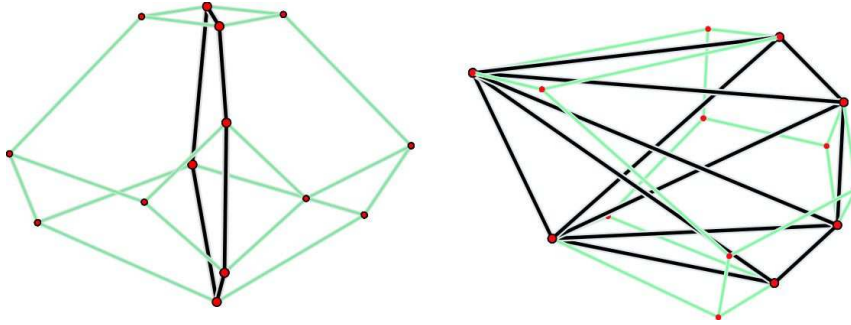}
         \end{center}
         \caption[]{Coordinates obtained from symmetric triangulations with symmetric~$\ADG$ yield a generalized
                    associahedron of type~$B_{n}$ as shown on the left (black edges). If~$\ADG$ is not symmetric,
                    the convex hull does not even yield a polytope of the correct dimension as shown on the right
                    (black edges).}
         \label{fig:convex_hulls}
      \end{minipage}
      \end{center}
\end{figure}

\begin{figure}
      \begin{center}
      \begin{minipage}{0.95\linewidth}
         \begin{center}
         \begin{overpic}
            [width=\linewidth]{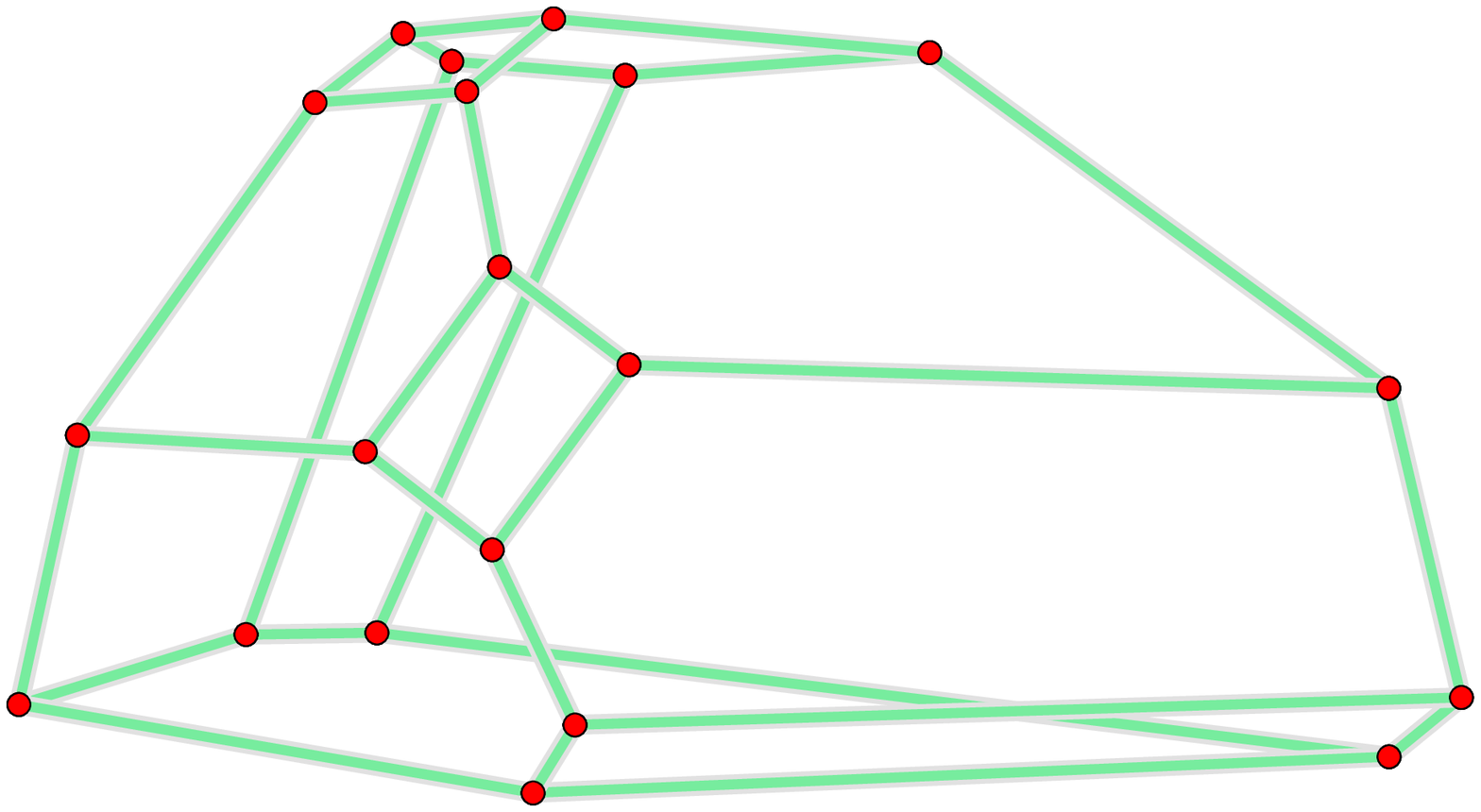}
            \put(46.5,68){\scriptsize $1$}
            \put(28,69){\scriptsize $2$}
            \put(44,54.5){\scriptsize $3$}
            \put(12,37){\scriptsize $4$}
            \put(35,34.5){\scriptsize $5$}
            \put(5,13){\scriptsize $6$}
            \put(50,29.5){\scriptsize $7$}
            \put(49.5,4.5){\scriptsize $8$}
            \put(50,14.5){\scriptsize $9$}
            \put(127,8){\scriptsize $10$}
            \put(134.5,16){\scriptsize $11$}
            \put(60,43){\scriptsize $12$}
            \put(128.5,44){\scriptsize $13$}
            \put(36,19){\scriptsize $14$}
            \put(24.5,19){\scriptsize $15$}
            \put(86.5,74.5){\scriptsize $16$}
            \put(59.5,69){\scriptsize $17$}
            \put(39.1,71.8){\scriptsize $18$}
            \put(35,76){\scriptsize $19$}
            \put(51,79){\scriptsize $20$}
         \end{overpic}
         \end{center}
         \[
           \begin{array}{ccc}
              \text{label} & \text{coordinate} & \text{triangulation}\\\hline
                     1            &   (3,5,6,1,2,4)   & \begin{minipage}{0cm}
                                                           \rule{-6mm}{12mm}\begin{overpic}
                                                                              [height=10mm]{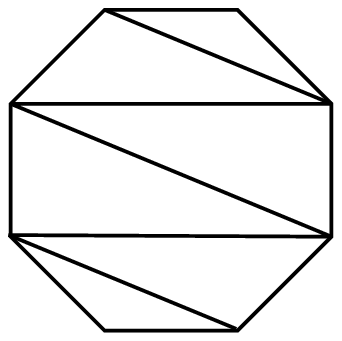}
                                                                              \put(1,6.8){\tiny $0$}
                                                                              \put(1,3){\tiny $1$}
                                                                              \put(3,9){\tiny $2$}
                                                                              \put(7,9){\tiny $3$}
                                                                              \put(3,1){\tiny $4$}
                                                                              \put(7,1){\tiny $5$}
                                                                              \put(9,3){\tiny $6$}
                                                                              \put(9,6.8){\tiny $7$}
                                                                           \end{overpic}
                                                        \end{minipage} \\
                     2            &   (3,6,5,2,1,4)   & \begin{minipage}{0cm}
                                                           \rule{-6mm}{12mm}\begin{overpic}
                                                                              [height=10mm]{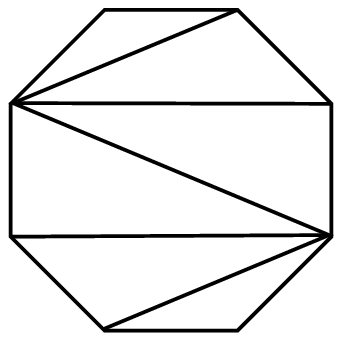}
                                                                              \put(1,6.8){\tiny $0$}
                                                                              \put(1,3){\tiny $1$}
                                                                              \put(3,9){\tiny $2$}
                                                                              \put(7,9){\tiny $3$}
                                                                              \put(3,1){\tiny $4$}
                                                                              \put(7,1){\tiny $5$}
                                                                              \put(9,3){\tiny $6$}
                                                                              \put(9,6.8){\tiny $7$}
                                                                           \end{overpic}
                                                        \end{minipage} \\
                     3            &   (2,4,6,1,3,5)   & \begin{minipage}{0cm}
                                                           \rule{-6mm}{12mm}\begin{overpic}
                                                                              [height=10mm]{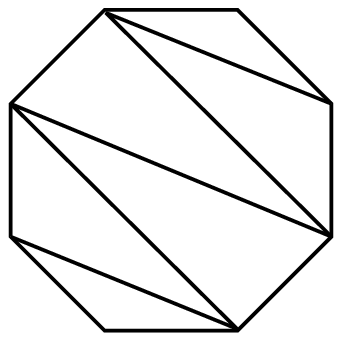}
                                                                              \put(1,6.8){\tiny $0$}
                                                                              \put(1,3){\tiny $1$}
                                                                              \put(3,9){\tiny $2$}
                                                                              \put(7,9){\tiny $3$}
                                                                              \put(3,1){\tiny $4$}
                                                                              \put(7,1){\tiny $5$}
                                                                              \put(9,3){\tiny $6$}
                                                                              \put(9,6.8){\tiny $7$}
                                                                           \end{overpic}
                                                        \end{minipage} \\
                     4            &   (1,6,3,4,1,6)   & \begin{minipage}{0cm}
                                                           \rule{-6mm}{12mm}\begin{overpic}
                                                                              [height=10mm]{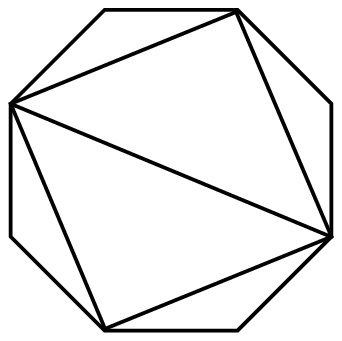}
                                                                              \put(1,6.8){\tiny $0$}
                                                                              \put(1,3){\tiny $1$}
                                                                              \put(3,9){\tiny $2$}
                                                                              \put(7,9){\tiny $3$}
                                                                              \put(3,1){\tiny $4$}
                                                                              \put(7,1){\tiny $5$}
                                                                              \put(9,3){\tiny $6$}
                                                                              \put(9,6.8){\tiny $7$}
                                                                           \end{overpic}
                                                        \end{minipage} \\
                     5            &   (1,4,5,2,3,6)   & \begin{minipage}{0cm}
                                                           \rule{-6mm}{12mm}\begin{overpic}
                                                                              [height=10mm]{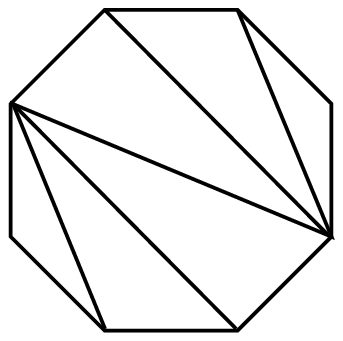}
                                                                              \put(1,6.8){\tiny $0$}
                                                                              \put(1,3){\tiny $1$}
                                                                              \put(3,9){\tiny $2$}
                                                                              \put(7,9){\tiny $3$}
                                                                              \put(3,1){\tiny $4$}
                                                                              \put(7,1){\tiny $5$}
                                                                              \put(9,3){\tiny $6$}
                                                                              \put(9,6.8){\tiny $7$}
                                                                           \end{overpic}
                                                        \end{minipage} \\
                     6            &   (1,6,-1,8,1,6)  & \begin{minipage}{0cm}
                                                           \rule{-6mm}{12mm}\begin{overpic}
                                                                              [height=10mm]{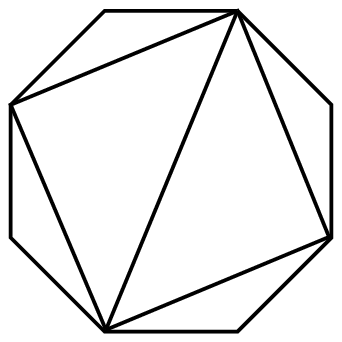}
                                                                              \put(1,6.8){\tiny $0$}
                                                                              \put(1,3){\tiny $1$}
                                                                              \put(3,9){\tiny $2$}
                                                                              \put(7,9){\tiny $3$}
                                                                              \put(3,1){\tiny $4$}
                                                                              \put(7,1){\tiny $5$}
                                                                              \put(9,3){\tiny $6$}
                                                                              \put(9,6.8){\tiny $7$}
                                                                           \end{overpic}
                                                        \end{minipage} \\
                     7            &   (1,3,5,2,4,6)   & \begin{minipage}{0cm}
                                                           \rule{-6mm}{12mm}\begin{overpic}
                                                                              [height=10mm]{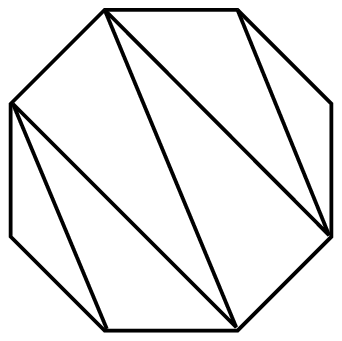}
                                                                              \put(1,6.8){\tiny $0$}
                                                                              \put(1,3){\tiny $1$}
                                                                              \put(3,9){\tiny $2$}
                                                                              \put(7,9){\tiny $3$}
                                                                              \put(3,1){\tiny $4$}
                                                                              \put(7,1){\tiny $5$}
                                                                              \put(9,3){\tiny $6$}
                                                                              \put(9,6.8){\tiny $7$}
                                                                           \end{overpic}
                                                        \end{minipage} \\
                     8            &   (1,2,3,4,5,6)   & \begin{minipage}{0cm}
                                                           \rule{-6mm}{12mm}\begin{overpic}
                                                                              [height=10mm]{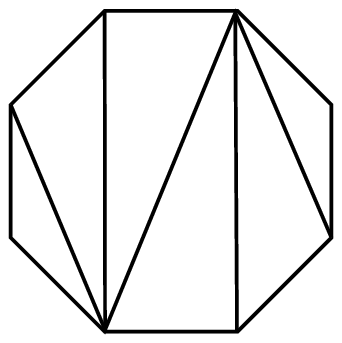}
                                                                              \put(1,6.8){\tiny $0$}
                                                                              \put(1,3){\tiny $1$}
                                                                              \put(3,9){\tiny $2$}
                                                                              \put(7,9){\tiny $3$}
                                                                              \put(3,1){\tiny $4$}
                                                                              \put(7,1){\tiny $5$}
                                                                              \put(9,3){\tiny $6$}
                                                                              \put(9,6.8){\tiny $7$}
                                                                           \end{overpic}
                                                        \end{minipage} \\
                     9            &   (1,2,4,3,5,6)   & \begin{minipage}{0cm}
                                                           \rule{-6mm}{12mm}\begin{overpic}
                                                                              [height=10mm]{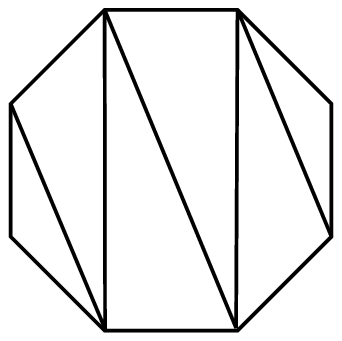}
                                                                              \put(1,6.8){\tiny $0$}
                                                                              \put(1,3){\tiny $1$}
                                                                              \put(3,9){\tiny $2$}
                                                                              \put(7,9){\tiny $3$}
                                                                              \put(3,1){\tiny $4$}
                                                                              \put(7,1){\tiny $5$}
                                                                              \put(9,3){\tiny $6$}
                                                                              \put(9,6.8){\tiny $7$}
                                                                           \end{overpic}
                                                        \end{minipage} \\
                     10           &   (6,-3,3,4,10,1) & \begin{minipage}{0cm}
                                                           \rule{-6mm}{12mm}\begin{overpic}
                                                                              [height=10mm]{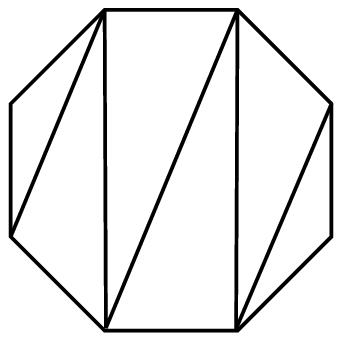}
                                                                              \put(1,6.8){\tiny $0$}
                                                                              \put(1,3){\tiny $1$}
                                                                              \put(3,9){\tiny $2$}
                                                                              \put(7,9){\tiny $3$}
                                                                              \put(3,1){\tiny $4$}
                                                                              \put(7,1){\tiny $5$}
                                                                              \put(9,3){\tiny $6$}
                                                                              \put(9,6.8){\tiny $7$}
                                                                           \end{overpic}
                                                        \end{minipage} \\
           \end{array}
           \qquad
           \begin{array}{ccc}
              \text{label} & \text{coordinate} & \text{triangulation}\\\hline
                     11           &   (6,-3,4,3,10,1) & \begin{minipage}{0cm}
                                                           \rule{-6mm}{12mm}\begin{overpic}
                                                                              [height=10mm]{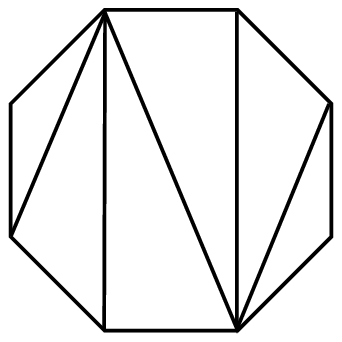}
                                                                              \put(1,6.8){\tiny $0$}
                                                                              \put(1,3){\tiny $1$}
                                                                              \put(3,9){\tiny $2$}
                                                                              \put(7,9){\tiny $3$}
                                                                              \put(3,1){\tiny $4$}
                                                                              \put(7,1){\tiny $5$}
                                                                              \put(9,3){\tiny $6$}
                                                                              \put(9,6.8){\tiny $7$}
                                                                           \end{overpic}
                                                        \end{minipage} \\
                     12           &   (2,3,6,1,4,5)   & \begin{minipage}{0cm}
                                                           \rule{-6mm}{12mm}\begin{overpic}
                                                                              [height=10mm]{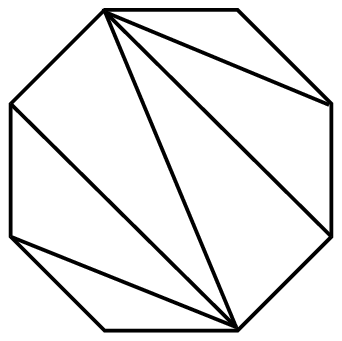}
                                                                              \put(1,6.8){\tiny $0$}
                                                                              \put(1,3){\tiny $1$}
                                                                              \put(3,9){\tiny $2$}
                                                                              \put(7,9){\tiny $3$}
                                                                              \put(3,1){\tiny $4$}
                                                                              \put(7,1){\tiny $5$}
                                                                              \put(9,3){\tiny $6$}
                                                                              \put(9,6.8){\tiny $7$}
                                                                           \end{overpic}
                                                        \end{minipage} \\
                     13           &   (6,-1,6,1,8,1)  & \begin{minipage}{0cm}
                                                           \rule{-6mm}{12mm}\begin{overpic}
                                                                              [height=10mm]{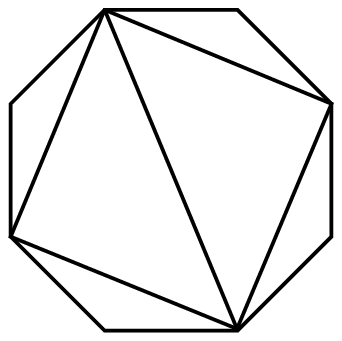}
                                                                              \put(1,6.8){\tiny $0$}
                                                                              \put(1,3){\tiny $1$}
                                                                              \put(3,9){\tiny $2$}
                                                                              \put(7,9){\tiny $3$}
                                                                              \put(3,1){\tiny $4$}
                                                                              \put(7,1){\tiny $5$}
                                                                              \put(9,3){\tiny $6$}
                                                                              \put(9,6.8){\tiny $7$}
                                                                           \end{overpic}
                                                        \end{minipage} \\
                     14           &   (6,5,-5,12,2,1) & \begin{minipage}{0cm}
                                                           \rule{-6mm}{12mm}\begin{overpic}
                                                                              [height=10mm]{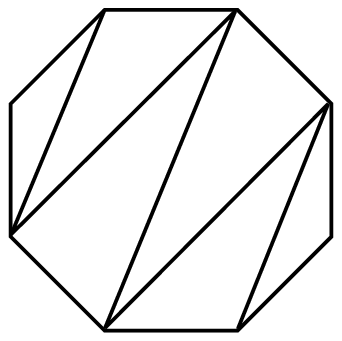}
                                                                              \put(1,6.8){\tiny $0$}
                                                                              \put(1,3){\tiny $1$}
                                                                              \put(3,9){\tiny $2$}
                                                                              \put(7,9){\tiny $3$}
                                                                              \put(3,1){\tiny $4$}
                                                                              \put(7,1){\tiny $5$}
                                                                              \put(9,3){\tiny $6$}
                                                                              \put(9,6.8){\tiny $7$}
                                                                           \end{overpic}
                                                        \end{minipage} \\
                     15           &   (5,6,-5,12,1,2) & \begin{minipage}{0cm}
                                                           \rule{-6mm}{12mm}\begin{overpic}
                                                                              [height=10mm]{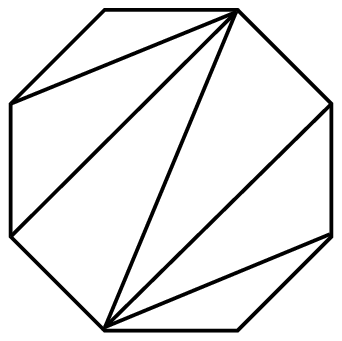}
                                                                              \put(1,6.8){\tiny $0$}
                                                                              \put(1,3){\tiny $1$}
                                                                              \put(3,9){\tiny $2$}
                                                                              \put(7,9){\tiny $3$}
                                                                              \put(3,1){\tiny $4$}
                                                                              \put(7,1){\tiny $5$}
                                                                              \put(9,3){\tiny $6$}
                                                                              \put(9,6.8){\tiny $7$}
                                                                           \end{overpic}
                                                        \end{minipage} \\
                     16           &   (6,3,6,1,4,1)   & \begin{minipage}{0cm}
                                                           \rule{-6mm}{12mm}\begin{overpic}
                                                                              [height=10mm]{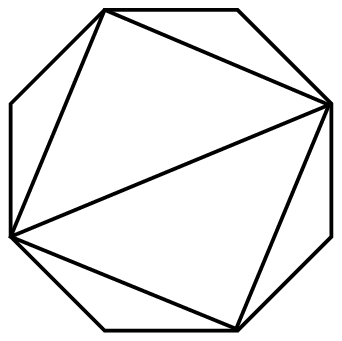}
                                                                              \put(1,6.8){\tiny $0$}
                                                                              \put(1,3){\tiny $1$}
                                                                              \put(3,9){\tiny $2$}
                                                                              \put(7,9){\tiny $3$}
                                                                              \put(3,1){\tiny $4$}
                                                                              \put(7,1){\tiny $5$}
                                                                              \put(9,3){\tiny $6$}
                                                                              \put(9,6.8){\tiny $7$}
                                                                           \end{overpic}
                                                        \end{minipage} \\
                     17           &   (6,5,4,3,2,1)   & \begin{minipage}{0cm}
                                                           \rule{-6mm}{12mm}\begin{overpic}
                                                                              [height=10mm]{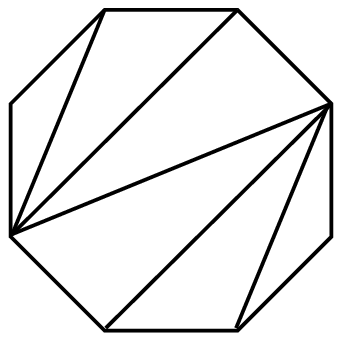}
                                                                              \put(1,6.8){\tiny $0$}
                                                                              \put(1,3){\tiny $1$}
                                                                              \put(3,9){\tiny $2$}
                                                                              \put(7,9){\tiny $3$}
                                                                              \put(3,1){\tiny $4$}
                                                                              \put(7,1){\tiny $5$}
                                                                              \put(9,3){\tiny $6$}
                                                                              \put(9,6.8){\tiny $7$}
                                                                           \end{overpic}
                                                        \end{minipage} \\
                     18           &   (5,6,4,3,1,2)   & \begin{minipage}{0cm}
                                                           \rule{-6mm}{12mm}\begin{overpic}
                                                                              [height=10mm]{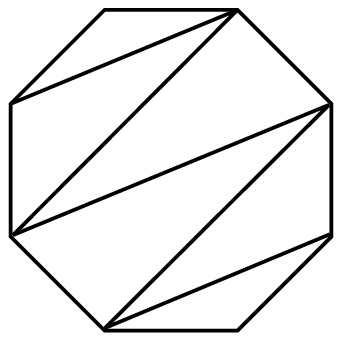}
                                                                              \put(1,6.8){\tiny $0$}
                                                                              \put(1,3){\tiny $1$}
                                                                              \put(3,9){\tiny $2$}
                                                                              \put(7,9){\tiny $3$}
                                                                              \put(3,1){\tiny $4$}
                                                                              \put(7,1){\tiny $5$}
                                                                              \put(9,3){\tiny $6$}
                                                                              \put(9,6.8){\tiny $7$}
                                                                           \end{overpic}
                                                        \end{minipage} \\
                     19           &   (4,6,5,2,1,3)   & \begin{minipage}{0cm}
                                                           \rule{-6mm}{12mm}\begin{overpic}
                                                                              [height=10mm]{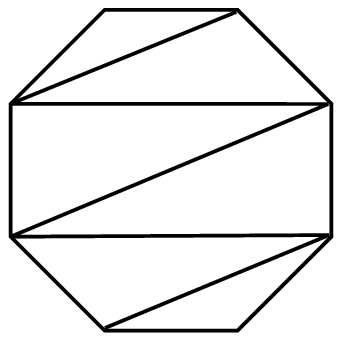}
                                                                              \put(1,6.8){\tiny $0$}
                                                                              \put(1,3){\tiny $1$}
                                                                              \put(3,9){\tiny $2$}
                                                                              \put(7,9){\tiny $3$}
                                                                              \put(3,1){\tiny $4$}
                                                                              \put(7,1){\tiny $5$}
                                                                              \put(9,3){\tiny $6$}
                                                                              \put(9,6.8){\tiny $7$}
                                                                           \end{overpic}
                                                        \end{minipage} \\
                     20           &   (4,5,6,1,2,3)   & \begin{minipage}{0cm}
                                                           \rule{-6mm}{12mm}\begin{overpic}
                                                                              [height=10mm]{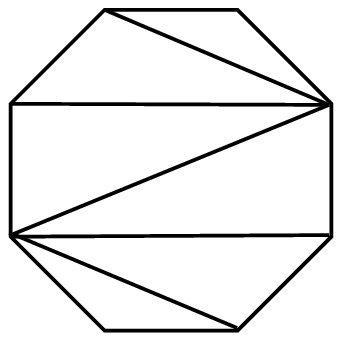}
                                                                              \put(1,6.8){\tiny $0$}
                                                                              \put(1,3){\tiny $1$}
                                                                              \put(3,9){\tiny $2$}
                                                                              \put(7,9){\tiny $3$}
                                                                              \put(3,1){\tiny $4$}
                                                                              \put(7,1){\tiny $5$}
                                                                              \put(9,3){\tiny $6$}
                                                                              \put(9,6.8){\tiny $7$}
                                                                           \end{overpic}
                                                        \end{minipage} \\
           \end{array}
         \]
         \caption[]{The vertex labels of this realization of the generalized associahedron of type~$B_{3}$
                    are decoded into coordinates and triangulations in table below. The corresponding
                    orientation of $\BCG_3$ is  obtained by directing the edges from left to right.}
         \label{fig:b3_associahedron}
      \end{minipage}
      \end{center}
\end{figure}

\noindent
A proof of Theorem~\ref{thm:AssB} is be given in Section~\ref{sec:RealizCyclo}, examples are shown in
Figures~\ref{fig:convex_hulls} and~\ref{fig:b3_associahedron}. The latter shows a realization 
of~$\Ass(\BCG_{3})$ together with a table of the coordinates of its vertices and the 
corresponding triangulations of the labelled octagon. We emphasize that Theorem~\ref{thm:AssB} is not 
true if the orientation~$\ADG$ is not symmetric as also visualized in Figure~\ref{fig:convex_hulls}: the 
obtained convex hull does not have the correct dimension, is not simple and  has triangular faces.

\medskip
The hyperoctahedral group~$W_{n}$ is the subgroup of~$S_{2n}$ that consists of all
permutations~$\sigma$ with the property $\sigma(2n+1-i)+\sigma(i)=2n+1$ for all
$i\in [n]$. As a Coxeter group of type $B_{n}$, the hyperoctahedral group is
generated by the simple transpositions $s_i:= \tau_i \tau_{2n-i}$, $i\in [n-1]$,
and the transposition~$t=\tau_{n}$. The Coxeter graph $\BCG_{n}$ is
\begin{figure}[h]
      \psfrag{t}{$t$}
      \psfrag{s1}{$s_{1}$}
      \psfrag{s2}{$s_{2}$}
      \psfrag{sn}{$s_{n-1}$}
      \psfrag{4}{$4$}
      \psfrag{dots}{$\ldots$}
      \begin{center}
            \includegraphics[width=8cm]{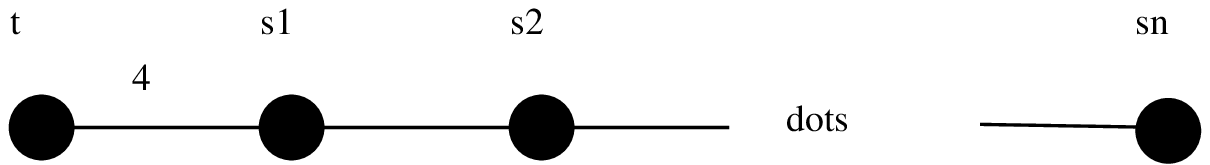}
      \end{center}
\end{figure}

\noindent
There is a bijection between the orientations of~$\BCG_{n}$ and the symmetric
orientations of~$\ACG_{2n-1}$. Let~$\BDG$ be an orientation of~$\BCG_{n}$,
then we construct an orientation of~$\ACG_{2n-1}$ by putting the
orientation~$\BDG$ on the subgraph of~$\ACG_{n-1}$ that consists
of the vertices $\tau_{n},\tau_{n+1},\dots,\tau_{2n-1}$, and by
completing the orientation symmetrically with respect to~$\tau_{n}$.
For convenience, we sometimes refer a symmetric orientation on~$\ACG_{2n-1}$
as~$\BDG$. For instance, the orientation
\begin{figure}[!h]
      \psfrag{t}{$t$}
      \psfrag{s1}{$s_{1}$}
      \psfrag{s2}{$s_{2}$}
      \psfrag{s3}{$s_{3}$}
      \psfrag{s4}{$s_{4}$}
      \psfrag{4}{$4$}
      \psfrag{dots}{$\ldots$}
      \begin{center}
            \includegraphics[width=8cm]{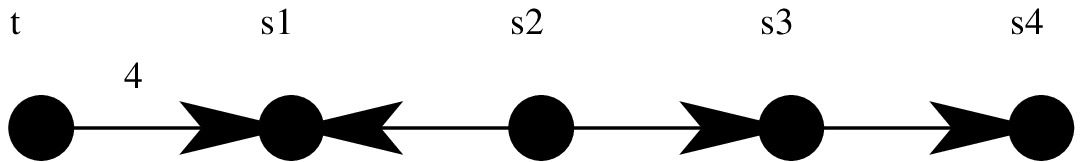}
      \end{center}
\end{figure}
$ $\\
\noindent
on $\BCG_5$ gives the following orientation on $\ACG_9$\\[-0.2cm]
\begin{figure}[!h]
      \psfrag{0}{$\tau_{1}$}
      \psfrag{1}{$\tau_{2}$}
      \psfrag{2}{$\tau_{3}$}
      \psfrag{3}{$\tau_{4}$}
      \psfrag{4}{$\tau_{5}$}
      \psfrag{5}{$\tau_{6}$}
      \psfrag{6}{$\tau_{7}$}
      \psfrag{7}{$\tau_{8}$}
      \psfrag{8}{$\tau_{9}$}
      \begin{center}
            \includegraphics[width=8cm]{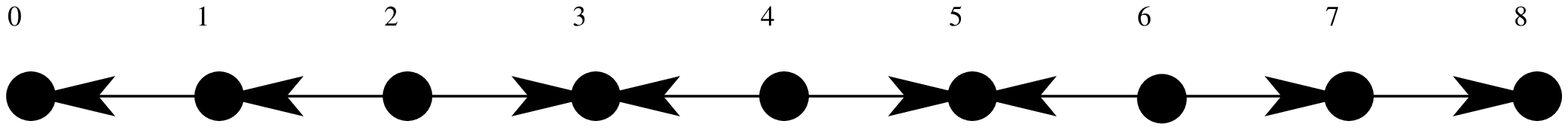}
      \end{center}
\end{figure}
\begin{figure}
      \begin{center}
      \begin{minipage}{0.95\linewidth}
         \begin{center}
            \includegraphics[height=8cm]{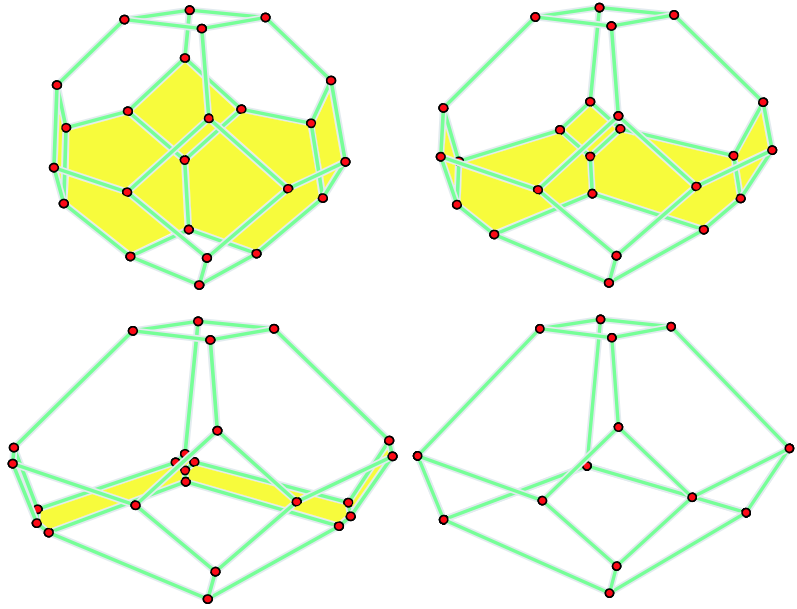}
         \end{center}
         \caption[]{The facets that correspond to non-admissible inequalities for the symmetric 
                   orientation~$\ADG$ of Figure~\ref{fig:a3_associahedra} (left associahedron) are
                   coloured. The four pictures show the process of removing these hyperplanes from the 
                   $A_{3}$-permutahedron (upper left) to the associahedron (bottom right).}
         \label{fig:removal_a_la_cambrian}
      \end{minipage}
      \end{center}
\end{figure}
\begin{figure}
      \begin{center}
      \begin{minipage}{0.95\linewidth}
         \begin{center}
            \includegraphics[height=8cm]{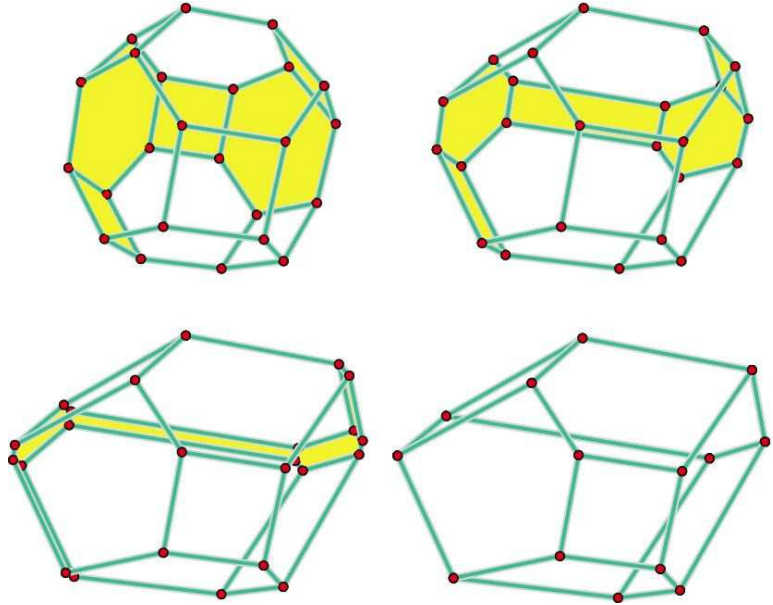}
         \end{center}
         \caption[]{The facets that correspond to non-admissible inequalities for the non-symmetric 
                    orientation~$\ADG$ of Figure~\ref{fig:a3_associahedra} (right associahedron) are
                    coloured (the perspective has changed by roughly 90 degrees with respect to the vertical
                    direction for a better visualization). The four pictures show the process of removing these
                    hyperplanes from the $A_{3}$-permutahedron (upper left) to the associahedron (bottom right).}
         \label{fig:removal_a_la_tamari}
      \end{minipage}
      \end{center}
\end{figure}

Let~$\BDG$ be an orientation of~$\BCG_{n}$, or equivalently
a symmetric orientation of~$\ACG_{2n-1}$. Denote by~$\Phi^B_\BDG$
the restriction~$\Phi_{\BDG | W_{n}}$ of the map~$\Phi_\BDG$ to~$W_{n}$. Then
\[
  \Phi_\BDG^B: W_n \longrightarrow \SOT_{2n+2}^B
\]
is surjective.

N.~Reading showed that~$\Phi^B_\BDG$  is a surjective lattice homomorphism
from the weak order lattice on~$W_{n}$ to a cambrian lattice of type~$B_{n}$.
Again, the undirected Hasse diagram of each cambrian lattice of type~$B_{n}$ is 
combinatorially equivalent to the $1$-skeleton of~$\Ass(\BCG_{n})$,~\cite[Theorem 1.3]{reading}. 
The permutahedron~$\Perm(\BCG_{n})$ of type~$B_{n}$ is the convex hull of the points
\[
  M(\sigma)=(\sigma(1),\sigma(2),\dots,\sigma(2n))\in \R^{2n},
  \qquad \forall \sigma\in W_{n}\subset S_{2n} .
\]
The next two propositions show that the realizations of the cyclohedron given in 
Theorem~\ref{thm:AssB} have similar properties as the ones of the associahedron
given in Theorem~\ref{thm:AssA}: They are obtained by removing certain inequalities
from the inequalities for~$\Perm(\BCG_n)$ and the common vertices~$\Ass(\BCG_{n})$
and~$\Perm(\BCG_{n})$ are characterized in many ways.

\begin{prop}\label{prop:AssBfromPermB}
   Fix an orientation~$\BDG$. The associahedron~$\Ass(\BCG_n)$ of Theorem~\ref{thm:AssB} 
   is given by a subset of the inequalities for the permutahedron~$\Perm(\BCG_{n-1})$. 
   These inequalities are determined by the image under~$K_{\BDG}$ of the diagonals of 
   the $(2n+2)$-gon labelled according to~$\BDG$.
\end{prop}

\begin{prop}\label{prop:AssBandPermB}
   Fix an orientation $\BDG$ on $\BCG_{n}$ and let $T\in \SOT_{2n+2}$ be centrally symmetric
   and $\sigma\in W_{n} \subset S_{2n}$. The following statements are equivalent:
   \begin{compactenum}[(a)]
      \item $M_{\BDG}(T) = M(\sigma)$,
      \item $\Phi^B_{\BDG}(\sigma)=T$ and the diagonals of~$T$ can be labelled such that
            \[
              \varnothing \subset K_{\ADG}(D_1) \subset \ldots \subset K_{\ADG}(D_{2n-1}) \subset [2n]
            \]
            is a sequence of strictly increasing nested sets.
      \item $(\Phi^B_{\BDG})^{-1}(T)=\{\sigma\}$,
      \item $\Phi^B_{\BDG}(\sigma) = T$ and for each~$i\in[2n]$ we have~$p_{\ell}^{T}(i)=1$ or~$p_{r}^{T}(i)=1$.
   \end{compactenum}
\end{prop}

The proofs are given in Section~\ref{sec:RealizCyclo}.

\subsection{Concerning the proofs}

The general idea to prove these results is to follow Loday's strategy: We start with 
a classical H-representation of~$\Perm(\ACG_{n-1})$, i.e. a representation by (in)equalities.
Then we identify among all defining inequalities the {\em $\ADG$-admissible} ones.
These are in bijection to the diagonals of the labelled $(n+2)$-gon and are precisely the
inequalities of an H-representation of~$\Ass(\ACG_{n-1})$. Finally, we show that
the intersection of all $\ADG$-admissible half spaces whose diagonals define a
triangulation $T\in \SOT_{n+2}$ is the point $M_{\ADG}(T)$. The process of removing
the non-admissible hyperplanes is visualized in Figures~\ref{fig:removal_a_la_cambrian}
and~\ref{fig:removal_a_la_tamari}. The facets supported by non-admissible inequalities
are shaded.

In his proof, Loday used two vital tools: A precise description of the admissible half spaces 
given by Stasheff,~\cite[Appendix]{stasheff2}, and the fact that any planar binary tree can
be cut into two planar binary trees. The latter piece of information gives rise to an inductive
argument.

In Section~\ref{sec:ProofAss}, we generalize Stasheff's H-representation
of~$\Ass(\ACG_{n-1})$ for all orientations of~$\ACG_{n-1}$, using results of
Reading,~\cite{reading}. But the induction of Loday does not generalize to our
set-up. We give a different proof that uses bistellar flips on triangulations
(i.e. flips of diagonals).

The permutahedron~$\Perm(\BCG_{n})$ of type~$B$ can be obtained by intersecting
the permutahedron~$\Perm(\ACG_{2n-1})$ with ``type-$B$-hyperplanes''. If the
orientation~$\ADG$ of~$\ACG_{2n-1}$ is symmetric, we conclude that the following
diagram is commutative:
$$
\xymatrix{
&\Ass(\ACG_{2n-1})\ar[dr]^\star\ar[dl]^\diamond&\\
\Perm(\ACG_{2n-1})\ar[dr]^\star&&\Ass(\BCG_{n})\ar[dl]^\diamond\\
&\Perm(\BCG_{n})&}
$$
The symbol $\diamond$ indicates that we intersect the starting polytope with all
non-admissible half spaces, and the symbol $\star$ indicates that we intersect the
starting polytope with the ``type-$B$-hyperplanes''. This gives the general idea
of the proof for type~$B$.

\section{H-representations of the associahedron and proofs for Subsection~\ref{sse:IntroAsso}}\label{sec:ProofAss}

\begin{figure}[t]
      \begin{center}
      \begin{minipage}{0.95\linewidth}
         \begin{center}
         \begin{overpic}
            [width=0.7\linewidth]{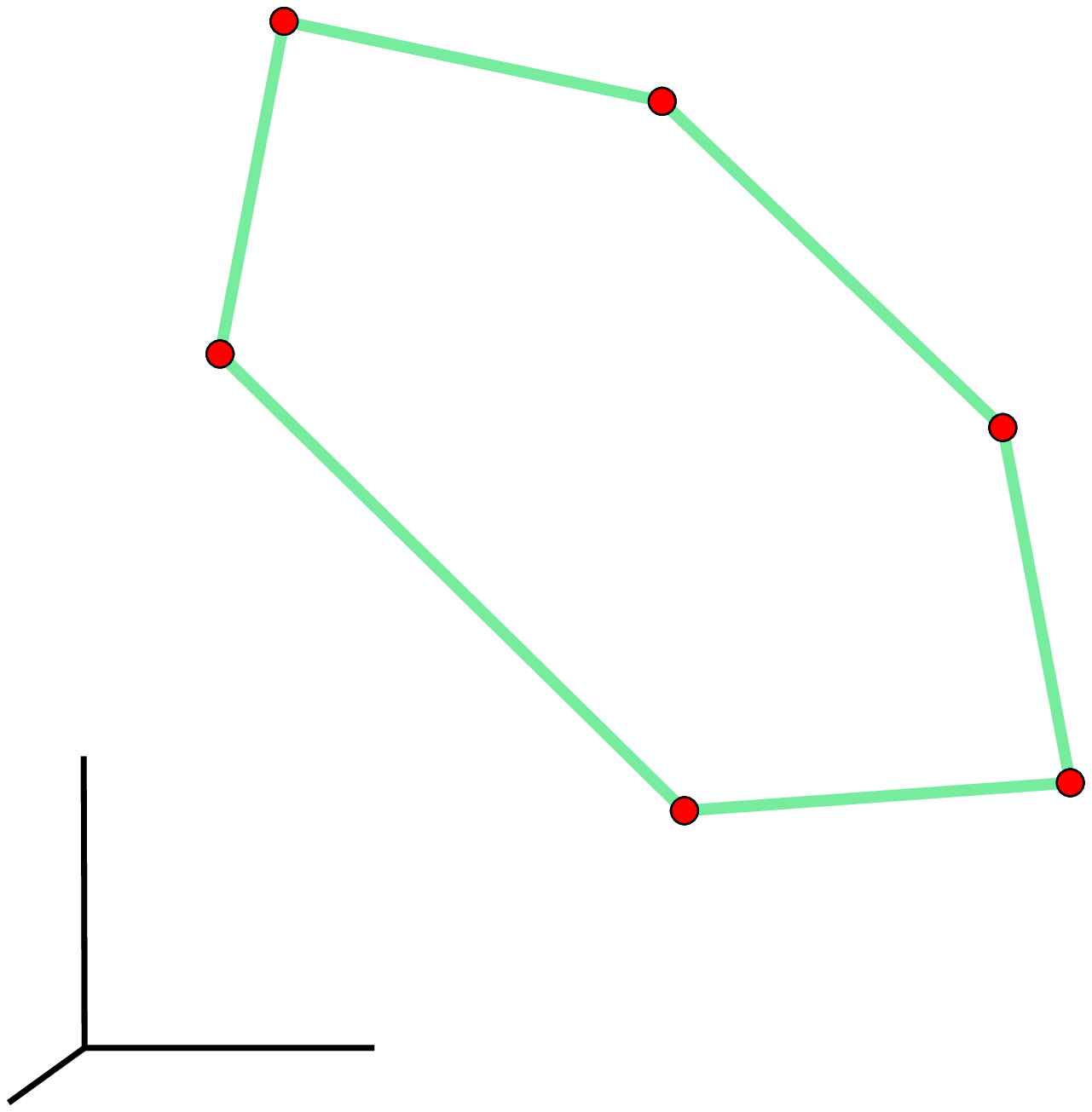}
            \put(4,47){$M(123)=(1,2,3)$}
            \put(8,67){$M(132)=(1,3,2)$}
            \put(59,64){$M(231)=(2,3,1)$}
            \put(79,45){$M(321)=(3,2,1)$}
            \put(84,21.5){$M(312)=(3,1,2)$}
            \put(31,20){$M(213)=(2,1,3)$}
            \put(36,55){$\{ 1 \}$}
            \put(70,25){$\{ 2 \}$}
            \put(63,50){$\{ 3 \}$}
            \put(47,36){$\{ 1,2 \}$}
            \put(44,61){$\{ 1,3 \}$}
            \put(72,32){$\{ 2,3 \}$}
            \put(40,8){\tiny{$x_{1}$}}
            \put(24.5,25){\tiny{$x_{2}$}}
            \put(22.5,4){\tiny{$x_{3}$}}
         \end{overpic}
         \end{center}
         \caption[]{The convex hull of~$\set{M(\sigma)}{\sigma\in S_3}$ yields a $2$-dimensional permutahedron
                    in~$\R^3$ contained in the affine hypeplane~$H$ with $x_1+x_2+x_3=6$. The 
                    intersections~$H \cap H_{K}$ for $\varnothing \subset K \subset [n]$ are the lines 
                    defined  by the edges of~$\Perm(\ACG_2)$. The edges are labelled by the set~$K$.}
         \label{fig:a2_permutahedron}
      \end{minipage}
      \end{center}
\end{figure}

We start with a classical  H-representation of the permutahedron~$\Perm(\ACG_{n-1})$ with vertex
set~$\set{M(\sigma)}{\sigma\in S_n}$, see Figure~\ref{fig:a2_permutahedron} for~$\Perm(\ACG_{2})$. 
Firstly, we consider the hyperplane
\[
  H = \set{x \in \mathbb R^{n}}{\sum_{i \in [n]}x_{i} = \tfrac{n(n+1)}{2}}.
\]
Secondly, each non-empty proper subset $K\subset [n]$ with~$k:=|K|$ defines the closed half space
\[
   \mathscr H_{K}
    := \set{ x \in \mathbb R^n}
           {(n-k)\sum_{i \in K}x_{i} - k\sum_{i\in [n]\setminus K}x_{i} + \tfrac{nk(n-k)}{2} \geq 0}.
\]
The open half space~$\mathscr H^{+}_{K}$ and the hyperplane~$H_{K}$
are defined by strict inequality and equality respectively. The negative half space~$\mathscr H^{-}_{K}$
is the complement of $\mathscr H_{K}$ in $\mathbb R^n$.
Now the permutahedron can be described as:
\[
  \Perm(\ACG_{n-1}) = H \cap \bigcap_{\varnothing \neq K \subset [n]}\mathscr H_{K}.
\]
Moreover, $M(\sigma)\in H_K$ if and only if $\sigma^{-1}(\big[|K|\big])=K$, see for
instance~\cite[\S2.2]{loday}. In other words,
\begin{equation}\label{equ:PointHyperplanePermA}
   \{M(\sigma) \} = H \cap
                    \bigcap_{{\scriptstyle \varnothing\neq K \subset [n] \atop \scriptstyle K=\sigma^{-1} ([|K|])}} H_K .
\end{equation}

Let~$P$ be the~$(n+2)$-gon labelled according to a given orientation~$\ADG$ of~$\ACG_{n-1}$.
We now describe an injective map $K_{\ADG}$ from the set of diagonals of~$P$ to the set of
non-empty proper subsets of~$[n]$. Set $\overline\Do_{\ADG}=\Do_{\ADG}\cup \{0,n+1\}$. For
a diagonal $D=\{a,b\}$, $0\leq a< b \leq n+1$, we define
\[
 K_{\ADG}(D) := \begin{cases}
                   \set{i \in \Do_{\ADG}}{a < i < b}
                                    &\textrm{if } a,b \in \overline\Do_{\ADG}\\
                   \set{i \in \Do_{\ADG}}{a < i} \cup \set{i\in \Up_{\ADG}}{b \leq i}
                                    &\textrm{if } a \in\overline\Do_{\ADG}, b \in \Up_{\ADG}\\
                   \Do_{\ADG} \cup \set{i\in \Up_{\ADG}}{i \leq a \text{ or } b \leq i}
                                    &\textrm{if } a, b \in \Up_{\ADG}\\
                   \set{i \in \Do_{\ADG}}{b > i} \cup \set{i \in \Up_{\ADG}}{a \geq i}
                                    &\textrm{if } a \in \Up_{\ADG}, b \in \overline\Do_{\ADG}.
                \end{cases}
\]
In other words, $K_\ADG(D)$ is the subset of $[n]$ obtained by
reading counterclockwise the labels of $P$ starting from $a$ and
ending with $b$, and by removing $0$, $n+1$, and $\{a,b\}\cap
\Do_\ADG$.

\begin{defn} \label{def:admissible_half_space}
   \textnormal{Fix an orientation~$\ADG$ of~$\ACG_{n-1}$. The half space~$\mathscr H_{K}$ is
           {\em $\ADG$-admissible} if there is a diagonal~$D$ of~$P$ such that $K=K_{\ADG}(D)$. }
\end{defn}
For instance, the $\ADG$-admissible half-spaces for the symmetric orientation~$\ADG$ corresponding 
to the realization on the left of Figure~\ref{fig:a3_associahedra} correspond to the subsets
$$
\{1\},\ \{3\},\ \{4\},\ \{1,2\},\ \{1,3\},\ \{3,4\},\ \{1,2,3\},\ \{1,3,4\},\ \{2,3,4\},
$$
while the admissible half spaces of the other realization in Figure~\ref{fig:a3_associahedra} 
correspond to
$$
 \{1\},\ \{2\},\ \{4\},\ \{1,2\},\ \{1,4\},\ \{2,3\},\ \{1,2,3\},\ \{1,2,4\},\ \{1,3,4\}.
$$

\medskip

We first proof a sequence of lemmas and corollaries to obtain in Proposition~\ref{prop:PropAss}
a better understanding of the relationship between the points~$M_{\ADG}(T)$ on the one hand 
and~$H$,~$H_{K}$, and~$\mathscr H^{+}_{K}$ on the other hand. The aim is to show that
\begin{compactenum}
   \item $M_{\ADG}(T)$ is contained in~$H$ for every triangulation~$T$;
   \item $D$ is a diagonal of a triangulation~$T$ if and only if $M_{\ADG}(T) \in H_{K_{\ADG}(D)}$;
   \item If~$D$ is not a diagonal of the triangulation~$T$ then $M_{\ADG}(T) \in \mathscr H^{+}_{K_{\ADG}(D)}$.
\end{compactenum}

\medskip
\noindent
We first show that certain half spaces are admissible for any orientation~$\ADG$ of~$\ACG_{n-1}$.

\begin{lem}\label{lem:top_bottom_admissible}
   For any orientation~$\ADG$ of~$\ACG_{n-1}$ the sets~$S_{u}=\{1,2, \ldots , u\}$ for $1 \leq u \leq n-1$
   and~$\widetilde S_{v}=\{ n,n-1, \ldots , n-v\}$ for $0 \leq v \leq n-2$ yield $\ADG$-admissible
   half spaces~$\mathscr H_{S_{u}}$ and~$\mathscr H_{\widetilde S_{v}}$.
\end{lem}

\begin{proof}  
   Denote the elements of $\Do_{\ADG}$ by $i_{1}=1 < i_{2} < \ldots < i_{\alpha}=n$ and the 
   elements of~$\Up_{\ADG}$ by $j_{1}< \ldots < j_{\beta}$. Let $1\leq u\leq n-1$. Let $j_k$ 
   be the greatest integer in $(S_u\cap\Up_\ADG)\cup\{0\}$ and $i_\ell$ be the greatest integer 
   in $(S_u\cap\Do_\ADG)$. Observe that $i_\ell < n$ since $i_\ell \in S_u$. Therefore $j_k<i_{l+1}$ 
   and the diagonal $\{j_k,i_{l+1}\}$ is mapped to $S_u$ under $K_\ADG$. Proceed similarly with $\widetilde S_v$.
\end{proof}

Both associahedra considered in Figure~\ref{fig:a3_associahedra} have a vertex~$(1,2,3,4)$ and $(4,3,2,1)$. The 
first vertex corresponds in both cases to the triangulation with diagonals $\{1\},\ \{1,2\},\ \{1,2,3\}$ and 
the second vertex corresponds in both cases to the triangulation with diagonals $\{4\},\ \{4,3\},\ \{4,3,2\}$. 
This is true in general. More precisely, we have the following corollary.

\begin{cor}\label{cor:top_bottom_vertex}
   For any orientation~$\ADG$ of~$\ACG_{n-1}$ there are triangulations~$T$ and~$\widetilde T$ of the labelled
   $(n+2)$-gon such that $M_{\ADG}(T)=(1,2,,\ldots,n)$ and $M_{\ADG}(\widetilde T)=(n,n-1,\ldots,1)$.
\end{cor}
\begin{proof}
   The diagonals described in the proof of Lemma~\ref{lem:top_bottom_admissible} to obtain the sets~$S_{u}$ yield 
   a triangulation~$T$ with $M_{\ADG}(T)=(1,2,\ldots,n)$ and the diagonals for the sets~$\widetilde S_{u}$ yield
   a triangulation~$\widetilde T$ with $M_{\ADG}(\widetilde T)=(n,n-1,\ldots,1)$.
\end{proof}

\begin{defn}
   A triangulation~$T\in \SOT_{n+2}$ \emph{refines} a given diagonal~$D$ if this diagonal~$D$
   is used in the triangulation~$T$. We write~$T\prec D$ in this situation.
\end{defn}

\begin{lem}\label{lem:LeAss1}
   Fix an orientation $\ADG$ on $\ACG_{n-1}$.
   \begin{compactenum}[1.]
      \item There is a triangulation $T\in \SOT_{n+2}$  with $M_{\ADG}(T) = (x_{1}, \ldots, x_{n})$
            such that
            \[
              \sum_{i \in [n]} x_{i} = \tfrac{n(n+1)}{2}.
            \]
      \item For each diagonal $D$ with $d := |K_{\ADG}(D)|$, there is a triangulation $T\in \SOT_{n+2}$
             refining $D$ with
            $M_{\ADG}(T) = (x_{1}, \ldots, x_{n})$ such that
            \[
              \sum_{i \in K_\ADG(D)} x_{i} = \tfrac{d(d+1)}{2}.
            \]
   \end{compactenum}
\end{lem}
\begin{proof}
   Denote the elements of~$\overline\Do_{\ADG}$ by
   $i_{0}=0<i_{1}=1 < i_{2} < \ldots < i_{\alpha-1} < i_{\alpha}=n < i_{\alpha+1}=n+1$ and the
   elements of $\Up_{\ADG}$ by $j_{1} <j_2< \ldots < j_{\beta}$. Observe that $\alpha+\beta=n$.
   \begin{compactenum}[1.]
      \item The triangulations~$T$ and~$\widetilde T$ of Corollary~\ref{cor:top_bottom_vertex} 
            have the desired property.
      \item We have to distinguish four cases. The aim is to produce a permutation~$\sigma\in S_n$ such
            that~$\Phi_\ADG(\sigma)$ refines a given diagonal~$D$. The first~$d$ elements of~$\sigma^{-1}$
            are precisely the elements of~$K_\ADG(D)$, therefore it is sufficient to specify the first~$d$
            elements of~$\sigma^{-1}$. This is what we shall do.
            \begin{compactenum}[\it i.]
               \item $D=\{ a,b\}$ with $a,b \in \overline\Do_{\ADG}$ and $a<b$.\\
                     Let~$u<v$ be such that~$i_{u}=a$ and~$i_{v}=b$. Then the desired triangulation is
                     obtained from any permutation~$\sigma$ when ~$\sigma^{-1}$ starts with the word
                     $i_{u+1}i_{u+2}\ldots i_{v-1}$. More precisely, we have~$x_{i_{u+1}}=1\cdot 1$,
                     $x_{i_{u+2}}= 2 \cdot 1$, $\ldots$, $x_{i_{v-1}}=(v-1-u) \cdot 1=d$, i.e.
                     $\sum_{i\in K_{\ADG}(D)}x_{i}=\tfrac{d(d+1)}{2}$.
               \item $D=\{ a,b\}$ with $a\in \overline\Do_{\ADG}$, $b\in \Up_{\ADG}$, and $a<b$.\\
                     Let $u$ and $v$ be such that $i_{u}=a$ and $j_{v}=b$. Then any permutation~$\sigma$
                     where~$\sigma^{-1}$ starts with $i_{\alpha}i_{\alpha-1}\ldots i_{u+1}j_{\beta}j_{\beta-1}\ldots j_{u}$
                     yields $\sum_{i\in K_{\ADG}(D)}x_{i}=\tfrac{d(d+1)}{2}$.
               \item $D=\{ a,b\}$ with $a\in \Up_{\ADG}$, $b\in \overline\Do_{\ADG}$, and $a<b$.\\
                     Let $u$ and $v$ be such that $i_{u}=b$ and $j_{v}=a$. Consider the coordinates obtained from
                     the triangulation associated to any permutation~$\sigma$ where $\sigma^{-1}$ starts with the word
                     $i_{1}i_{2}\ldots i_{u-1}j_{1}j_{2}\ldots j_{v}$.
               \item $D=\{ a,b\}$ with $a,b\in \Up_{\ADG}$ and $a<b$.\\
                     Let $u$ and $v$ be such that $j_{u}=a$ and $j_{v}=b$. Let~$\sigma$ be any permutation where~$\sigma^{-1}$
                     starts with the word~$i_{1}i_{2}\ldots i_{\alpha}j_{1}j_{2}\ldots j_{u}j_{\beta}j_{\beta-1}\ldots j_{v}$.
            \end{compactenum}
   \end{compactenum}
\end{proof}

Let $T\in \SOT_{n+2}$ a triangulation of the labelled $(n+2)$-gon~$P$ and~$\{a,c\}$ be a diagonal
of~$T$. There are two unique labels~$b$,~$d$ of~$T$ such that~$\{a,c\}$ is a diagonal of the quadrilateral
given by the edges $\{a,b\}$, $\{b,c\}$, $\{c,d\}$, and $\{a,d\}$ of~$T$. Hence the diagonal~$\{b,d\}$ is 
not an edge of~$T$. The {\em bistellar flip} of the diagonal~$\{a,c\}$ is the transformation which map~$T$ 
to~$T^{\prime}$ where $T^{\prime} \in \SOT_{n+2}$ is the triangulation obtained by replacing the diagonal~$\{a,c\}$
by the diagonal~$\{b,d\}$ in~$T$. For two triangulations $T, T^{\prime} \in \SOT_{n+2}$,
we write $T\dot\sim T^{\prime}$ if~$T^{\prime}$ can be obtained from~$T$ by a bistellar flip of
a diagonal of~$T$. The relation~$\dot\sim$ is symmetric. Denote by~$\sim$ the transitive and reflexive
closure of~$\dot\sim$. For any $T, T^{\prime} \in \SOT_{n+2}$, there is a sequence
$T=T_1,T_2,\dots,T_p=T^{\prime}$ of triangulations in $\SOT_{n+2}$ such that $T_i \dot\sim T_{i+1}$ for
all $i\in [p-1]$.

\begin{lem}\label{lem:LeAss1.2}
   Fix an orientation~$\ADG$ on~$\ACG_{n-1}$. Let~$T \in \SOT_{n+2}$ and~$D$ a diagonal of~$T$. Consider the
   triangulation~$T^{\prime}$ that is obtained from~$T$ by a bistellar flip from~$D$ to~$D^{\prime}$.
   Set $M_{\ADG}(T) = (x_{1}, \ldots, x_{n})$ and $M_{\ADG}(T^{\prime}) = (y_{1}, \ldots, y_{n})$.
  The vertices of the quadrilateral with diagonals~$D$ and~$D^{\prime}$ are labelled $a < b < c <d$. Then
  $x_{i}=y_{i}$ for all $i \in [n] \setminus \{b,c\}$ and   $x_{b} + x_{c} = y_{b} + y_{c}$.
\end{lem}
\begin{proof}
  It follows immediately from the definitions that~$x_{i}=y_{i}$ for all $i \in [n] \setminus \{b,c\}$.
  We have to show that $x_{b} + x_{c} = y_{b} + y_{c}$.
            There are~$4$ cases to distinguish: $b$ and $c$ are elements of~$\overline\Do_{\ADG}$ or~$\Up_{\ADG}$.
            \begin{compactenum}[i.]
               \item $b,c \in \overline\Do_{\ADG}$.\\
                     We have~$\mu_{c}(a) = \mu_{b}(a) + \mu_{c}(b)$,~$\mu_{b}(d) = \mu_{b}(c) + \mu_{c}(d)$, and
                     \begin{align*}
                        \mu_{b}(a)\mu_{b}(c) + \mu_{c}(a)\mu_{c}(d)
                               &= \mu_{b}(a)\mu_{b}(c) + [\mu_{b}(a)+\mu_{c}(b)]\mu_{c}(d)\\
                               &= \mu_{b}(a)[\mu_{b}(c) + \mu_{c}(d)] + \mu_{c}(b)\mu_{c}(d)\\
                               &= \mu_{b}(a)\mu_{b}(d) + \mu_{c}(b)\mu_{c}(d).
                     \end{align*}
                     If~$D = \{ a,c \}$ and~$D^{\prime} = \{ b,d \}$ we
                     have~$x_{b}+x_{c} = \mu_{b}(a)\mu_{b}(c) + \mu_{c}(a)\mu_{c}(d)$
                     and~$y_{b}+y_{c} = \mu_{b}(a)\mu_{b}(d) + \mu_{c}(b)\mu_{c}(d)$.
                     If~$D = \{ b,d \}$ and~$D^{\prime} = \{ a,c \}$ we
                     have~$y_{b}+y_{c} = \mu_{b}(a)\mu_{b}(c) + \mu_{c}(a)\mu_{c}(d)$
                     and~$x_{b}+x_{c} = \mu_{b}(a)\mu_{b}(d) + \mu_{c}(b)\mu_{c}(d)$.
               \item $b \in \overline\Do_{\ADG}$ and $c \in \Up_{\ADG}$.\\
                     We have~$\mu_{c}(a) = \mu_{c}(b) - \mu_{b}(a)$,~$\mu_{b}(c) = \mu_{b}(d) + \mu_{c}(d)$, and
                     \begin{align*}
                        \mu_{b}(a)\mu_{b}(d) + n + 1 - \mu_{c}(a)\mu_{c}(d)
                               &= \mu_{b}(a)\mu_{b}(d) + n + 1 - [\mu_{c}(b) - \mu_{b}(a)]\mu_{c}(d)\\
                               &= \mu_{b}(a)[\mu_{b}(d) + \mu_{c}(d)] + n + 1 - \mu_{c}(b)\mu_{c}(d)\\
                               &= \mu_{b}(a)\mu_{b}(c) + n + 1 - \mu_{c}(b)\mu_{c}(d)
                     \end{align*}
                     We have either~$D= \{ a,d \}$ and~$D^{\prime} = \{ b,c \}$
                     or $D= \{ b,c \}$ and $D^{\prime} = \{ a,d \}$.
                     Both cases imply~$x_{b}+x_{c} = y_{b}+y_{c}$.
               \item $b \in \Up_{\ADG}$ and $c \in \overline\Do_{\ADG}$.\\
                     We have~$\mu_{c}(a) = \mu_{c}(b) - \mu_{b}(a)$,~$\mu_{b}(c) = \mu_{b}(d) + \mu_{c}(d)$, and
                     \begin{align*}
                        n + 1 - \mu_{b}(a)\mu_{b}(d) + \mu_{c}(a)\mu_{c}(d)
                               &= n + 1 - \mu_{b}(a)\mu_{b}(d) + [\mu_{c}(b) - \mu_{b}(a))\mu_{c}(d)\\
                               &= n + 1 - \mu_{b}(a)[\mu_{b}(d) + \mu_{c}(d)] + \mu_{c}(b)\mu_{c}(d)\\
                               &= n + 1 - \mu_{b}(a)\mu_{b}(c) + \mu_{c}(b)\mu_{c}(d)
                     \end{align*}
                     We have either~$D= \{ a,d \}$ and~$D^{\prime} = \{ b,c \}$
                     or~$D= \{ b,c \}$ and~$D^{\prime} = \{ a,d \}$.
                     Both cases imply $x_{b}+x_{c} = y_{b}+y_{c}$.
               \item $b,c \in \Up_{\ADG}$.\\
                     We have~$\mu_{c}(a) = \mu_{c}(b) + \mu_{b}(a)$,~$\mu_{b}(c) = \mu_{b}(d) + \mu_{c}(d)$, and
                     \begin{align*}
                        n + 1 - \mu_{b}(a)\mu_{b}(c) + n + 1 &- \mu_{c}(a)\mu_{c}(d)\\
                               &= n + 1 - \mu_{b}(a)\mu_{b}(c) + n + 1 - [\mu_{c}(b) + \mu_{b}(a)]\mu_{c}(d)\\
                               &= n + 1 - \mu_{b}(a)[\mu_{b}(c) + \mu_{c}(d)] + n + 1 - \mu_{c}(b)\mu_{c}(d)\\
                               &= n + 1 - \mu_{b}(a)\mu_{b}(d) + n + 1 -\mu_{c}(b)\mu_{c}(d)
                     \end{align*}
                     We have either~$D= \{ a,c \}$ and~$D^{\prime} = \{ b,d \}$
                     or~$D= \{ b,d \}$ and~$D^{\prime} = \{ a,c \}$.
                     Both cases imply $x_{b}+x_{c} = y_{b}+y_{c}$.
            \end{compactenum}
\end{proof}

\begin{cor}\label{lem:LeAss2}
   Fix an orientation~$\ADG$ on~$\ACG_{n-1}$. Let $T \in \SOT_{n+2}$ and write $M_{\ADG}(T) = (x_{1}, \ldots, x_{n})$.
   \begin{compactenum}[1.]
      \item $\sum_{i \in [n]} x_{i}$ is invariant under bistellar flips of diagonals.
      \item Let~$D$ and~$D^{\prime}$ be distinct diagonals of~$T$, i.e.~$T$ refines both~$D$ and~$D^{\prime}$.
            Denote the triangulation obtained from a bistellar flip of~$D^{\prime}$ by~$T^{\prime}$ and
            $M_{\ADG}(T^{\prime})=(y_{1},\ldots , y_{n})$. Then
            \[
              \sum_{i \in K_{\ADG}(D)} y_{i} = \sum_{i \in K_\ADG(D)} x_{i}.
            \]
   \end{compactenum}
\end{cor}
\begin{proof}
   \begin{compactenum}[1.]
      \item Follows immediately from Lemma~\ref{lem:LeAss1.2}.
      \item The claim follows immediately from the first statement of this lemma: Let $a<b<c<d$ be the labels that
            define the quadrilateral for the bistellar flip of~$D^{\prime}$. Since~$T$ refines~$D$ and~$D^{\prime}$,
            we conclude that either $b,c \in K_{\ADG}(D)$ or $b,c \not\in K_{\ADG}(D)$.
   \end{compactenum}
\end{proof}

\noindent
A careful analysis of the proof of Lemma~\ref{lem:LeAss1.2} yields the following result.
\begin{cor}\label{cor:flip_inequality}
   Fix an orientation~$\ADG$ on~$\ACG_{n-1}$. Let~$T \in \SOT_{n+2}$ and~$D$ a diagonal of~$T$. Consider the
   triangulation~$T^{\prime}$ that is obtained from~$T$ by a bistellar flip from~$D$ to~$D^{\prime}$.
   Set $d = | K_{\ADG}(D) |$, $M_{\ADG}(T) = (x_{1}, \ldots, x_{n})$,
   and $M_{\ADG}(T^{\prime}) = (y_{1}, \ldots, y_{n})$. Then
   \[
     \sum_{i \in K_{\ADG}(D)} y_{i} > \sum_{i \in K_{\ADG}(D)} x_{i} = \tfrac{d(d+1)}{2}.
   \]
\end{cor}
\begin{proof}
   Again, we have to consider the quadrilateral spanned by the diagonals~$D$ and~$D^{\prime}$, its vertices are
   without loss of generality $a < b < c < d$. We only show the first case $b,c \in \overline\Do_{\ADG}$. The
   other cases are handled analogously.
   Suppose we flip from~$\{ a,c \}$ to~$\{ b,d \}$. Then $b \in K_{\ADG}(D)$ and $c \not \in K_{\ADG}(D)$.
   The claim follows from~$x_{b} < y_{b}$ as shown in the proof of Lemma~\ref{lem:LeAss1.2} since $\mu_b(c)<\mu_b(d)$.
\end{proof}

\begin{lem}\label{lem:LeAss3}
   Fix an orientation~$\ADG$ on $\ACG_{n-1}$. Let $T\in \SOT_{n+2}$ and write
    $M_{\ADG}(T) = (x_{1}, \ldots, x_{n})$. If~$T$ does not
    refine a given diagonal~$D$ with $d := |K_{\ADG}(D)|$ then
   \[
     \sum_{i \in K_\ADG(D)} x_{i} > \tfrac{d(d+1)}{2}.
   \]
\end{lem}
\begin{proof}
   Let~$u$ and~$v$ be the endpoints of~$D$ such that $u < v$.
   Since~$T$ does not refine~$D$, we have diagonals of~$T$ that intersect the line segment between~$u$ and~$v$
   in its relative interior. Let~$D_{1}, \ldots , D_{t}$ be all these diagonals ordered such
   that~$D \cap D_{i+1}$ (as intersection of line segments not as intersection of subsets of~$\{ 0,\ldots n+1 \}$)
   is closer to~$v$ than~$D \cap D_{i}$ for all~$i \in [t-1]$. Let~$u_{i}$
   and~$v_{i}$ denote the endpoints of~$D_{i}$ where~$u_{i} \in K_{\ADG}(D)$ and~$v_{i} \not \in K_{\ADG}(D)$
   for each $i \in [t]$.

   The strategy is now to flip the diagonal~$D_{1}$, then~$D_{2}, \ldots D_{t}$ to obtain a triangulation~$T^{\prime}$
   that refines~$D$. We show by induction on $t$ that $\sum_{i \in K_\ADG(D)} x_{i}$ decreases with each flip.

   We first remark that the special case~$t=1$ is covered by Corollary~\ref{cor:flip_inequality}
   (obtain a triangulation~$T^{\prime}$ that refines~$D$ by a bistellar flip from~$D_{1}$ to $D$).

   Suppose the claim is true for all $\bar t < t$.
   Apply a bistellar flip to the diagonal~$D_{1}$ of the quadrilateral~$\{a<b<c<d\}=\{ u, u_{1}, v_{1}, u_{2} , v_{2} \}$
    to obtain the triangulation~$T^{\prime}$
   with~$M_{\ADG}(T^{\prime}) = (y_{1}, \ldots , y_{n})$ and new diagonal~$D_{1}^{\prime}$
  (this is in fact a quadrilateral since there is no other diagonal intersecting $D$
    between $D_1$ and $D_2$ i.e. $u_{1}=u_{2}$ or $v_{1}=v_{2}$ ). In $T^\prime$, only
    $D_2,\ldots, D_t$ intersect the line between $u$ and $v$. We have
   $\sum_{i \in K_{\ADG}(D)}y_{i} > \tfrac{d(d+1)}{2}$ by induction, so it suffices to show
   $\sum_{i \in K_{\ADG}(D)}x_{i} \geq \sum_{i \in K_{\ADG}(D)}y_{i}$.

   From $D^{\prime}_{1} \cap D = \{ u \}$ we conclude that one of the following statements is true:
   \begin{compactenum}
      \item $K_{\ADG}(D^{\prime}_{1}) \subset K_{\ADG}(D)$,
      \item $K_{\ADG}(D^{\prime}_{1}) \supset K_{\ADG}(D)$,
      \item $K_{\ADG}(D^{\prime}_{1}) \cap    K_{\ADG}(D) =
      \varnothing$,
      \item $u=c\in \Up_\ADG$.
   \end{compactenum}

   Observe first that Corollary~\ref{cor:flip_inequality} implies that
   $\sum_{i \in K_{\ADG}(D_1^\prime)}x_{i} > \sum_{i \in K_{\ADG}(D_1^\prime)}y_{i}$.

   The first case implies that at least one of $b$ and $c$ is contained in $K_{\ADG}(D)$ (possibly both). From
   Lemma~\ref{lem:LeAss1.2}  we conclude $\sum_{i \in K_{\ADG}(D)}x_{i} \geq \sum_{i \in K_{\ADG}(D)}y_{i}$.

   The second case implies that either none, one, or both of $b,c$ are contained in $K_{\ADG}(D)$. If none or both are
   contained in $K_{\ADG}(D)$, we have $\sum_{i \in K_{\ADG}(D)}x_{i} = \sum_{i \in K_{\ADG}(D)}y_{i}$.
   If one of $b,c$ is contained in $K_{\ADG}(D)$, we have
   $\sum_{i \in K_{\ADG}(D)}x_{i} > \sum_{i \in K_{\ADG}(D)}y_{i}$ by  Lemma~\ref{lem:LeAss1.2}.

   The third case implies that $c=u$ and $u \in \overline \Do_{\ADG}$, i.e. $b,c \not \in K_{\ADG}(D)$ and. Hence
   we conclude $\sum_{i \in K_{\ADG}(D)}x_{i} = \sum_{i \in K_{\ADG}(D)}y_{i}$ by Lemma~\ref{lem:LeAss1.2}.

   The fourth case implies that $b,c$ are contained in
   $K_{\ADG}(D)$, then we have $\sum_{i \in K_{\ADG}(D)}x_{i} = \sum_{i \in K_{\ADG}(D)}y_{i}$ by Lemma~\ref{lem:LeAss1.2} again.
\end{proof}

As an consequence we obtain the following result:

\begin{prop}\label{prop:PropAss}
   Fix an orientation~$\ADG$ on~$\ACG_{n-1}$ and let $T\in \SOT_{n+2}$ and let $D$ be a diagonal. Then
   \begin{compactenum}
      \item $M_{\ADG}(T) \in H$,
      \item $T\prec D$ if and only if $M_{\ADG}(T) \in H_{K_{\ADG}(D)}$,
      \item $M_{\ADG}(T)\in \mathscr H^{+}_{K_{\ADG}(D)}$ if $T$ does not refine~$D$.
   \end{compactenum}
\end{prop}
\begin{proof}
   It is a well-known fact that any triangulation of a polygon can be transformed into any other
   triangulation by a sequence of bistellar flips. If both triangulations have a common diagonal,
   this sequence can be chosen in such that this diagonal is common to all intermediate triangulations.
   These remarks combined with Lemma~\ref{lem:LeAss1} and Corollary~\ref{lem:LeAss2} settle the first two statements.

   If $T$ does not refine~$D$, then write $M_{\ADG}(T)=(x_1,\dots,x_n)$ and $d:= |K_{\ADG}(D)|$. As $M_\ADG(T)\in H$,
   \[
     (n-d)\sum_{i \in K_\ADG(D)}x_{i} - d\sum_{i\in [n]\setminus K}x_{i} + \tfrac{nd(n-d)}{2}
     = n \sum_{i \in K_\ADG(D)}x_{i} - \tfrac{nd(d+1)}{2}
     > 0
   \]
   by Lemma~\ref{lem:LeAss3}. In other words, $M_{\ADG}(T)\in \mathscr H^{+}_{K_{\ADG}(D)}$.
\end{proof}

\begin{cor}\label{cor:AssAPoints}
   Fix an orientation~$\ADG$ on~$\ACG_{n-1}$ and let $T\in \SOT_{n+2}$. Then
   \[
     \left\{ M_\ADG(T) \right\} = H \cap \bigcap_{D\succ T} H_{K_{\ADG}(D)}.
   \]
\end{cor}
\begin{proof}
   It is clear that $\dim \left( H\cap \bigcap_{D\succ T} H_{K_{\ADG}(D)} \right) \leq 0$.
   But this intersection contains $M_\ADG(T)$ by Proposition~\ref{prop:PropAss}.
\end{proof}

\begin{thm}\label{cor:AssABound}
   The intersection of all $\ADG$-admissible half-spaces with~$H$ yields an associahedron
   with vertex set $\{M_\ADG(T),\, T\in\SOT_{n+2}\}$.
\end{thm}
\begin{proof}
   We first observe that the intersection of all admissible half spaces defines a bounded
   set in~$\R^{n}$. This follows  from the following facts:
   \begin{compactenum}
      \item From Lemma~\ref{lem:top_bottom_admissible},Corollary~\ref{cor:top_bottom_vertex},
            Proposition~\ref{prop:PropAss} and the H-representation of $\Perm(\ACG_{n-1})$, we conclude  
            that all half spaces~$\mathscr H_{K}$ that contain $(1,2,\ldots, n)$ or $(n,n-1,\ldots,1)$
            on their boundary~$H_{K}$ are admissible. The half spaces~$\mathscr H_{K}$ that contain 
            $(1,2,\ldots, n)$ on their boundary intersect with~$H$ in a cone~$C$ with apex
            $(1,2,,\ldots,n)$ of dimension~$\dim H$. Similarly, the half spaces~$\mathscr H_{K}$ that 
            contain $(n,n-1,\ldots,1)$ in their boundary intersect with~$H$ in a cone~$\widetilde C$ 
            with apex~$(n,n-1,\ldots,1)$ of dimension~$\dim H$. Since all these half spaces can be 
            partitioned into pairs $\mathscr H_{K}$ and $\mathscr H_{[n]\setminus K}$ where~$H_{K}$ 
            is parallel to $H_{[n]\setminus K}$ and $\mathscr H_{K} \supset \mathscr H^{-}_{[n]\setminus K}$,
            we conclude that the intersection $C \cap \widetilde C$ is a convex polytope.
      \item We intersect~$C \cap \widetilde C$ with all remaining admissible half spaces to obtain a convex
            polytope~$Q$ that contains~$\Perm(\ACG_{n-1})$.
   \end{compactenum}
   By Proposition~\ref{prop:PropAss}, we know that the
   vertex set~$V(Q)$ contains the set $\set{M_\ADG(T)}{T\in\SOT_{n+2}}$ and each vertex in
   $\set{M_\ADG(T)}{ T\in\SOT_{n+2}}$ is simple: it is contained in precisely
   $(n+2)-3 = n-1$ facet defining hyperplanes and in the interior of all other admissible half spaces. In particular,
   we conclude that each vertex of $\set{M_\ADG(T)}{T\in\SOT_{n+2}}$ is connected to precisely $(n-1)$ vertices
   of $\set{M_\ADG(T)}{T\in\SOT_{n+2}}$ by an edge: replace a defining hyperplane~$H_{1}$ of~$M_{\ADG}(T)$ by the
   hyperplane~$H_{2}$ that corresponds to the diagonal obtained from ``flipping~$H_{1}$ in~$T$''. This
   implies that all vertices of~$Q$ are contained in~$\set{M_\ADG(T)}{T\in\SOT_{n+2}}$ since the $1$-skeleton
   of a polytope is connected.

   Thus $Q$ is a simple polytope and its $1$-skeleton is the flip graph of an $(n+2)$-gon. This implies that~$Q$
   is an associahedron.
\end{proof}

\subsection{Proof of Theorem~\ref{thm:AssA}, Proposition~\ref{prop:AssAfromPermA}, and Proposition~\ref{prop:AssAandPermA}}
$ $

\medskip
\noindent
Theorem~\ref{thm:AssA} and Proposition~\ref{prop:AssAfromPermA} are immediate consequences of
Theorem~\ref{cor:AssABound}.

\begin{proof}[Proof of Proposition~\ref{prop:AssAandPermA}]
   $ $\\
   $(a)\Rightarrow (b)$: Denote the diagonals of~$T$ by $D_{1}, \ldots , D_{n-1}$. From 
      Equation~\ref{equ:PointHyperplanePermA}, Statement~$(a)$, and Corollary~\ref{cor:AssAPoints} we have
      \[
        H \cap \bigcap_{{\scriptstyle \varnothing\neq K \subset [n] \atop \scriptstyle K=\sigma^{-1} ([|K|])}} H_K
        \quad=\quad \left\{M(\sigma) \right\}
        \quad=\quad \left\{ M_\ADG(T) \right\} 
        \quad=\quad H \cap \bigcap_{D\succ T} H_{K_{\ADG}(D)}.
      \]
      Since $M(\sigma) \not\in H_{K}$ if $K \subset [n]$ and not of the type $\sigma^{-1}([r])$, $1 \leq r \leq n-1$,
      we may assume that 
      $K_{\ADG}(D_{i}) = \sigma^{-1}([i])$, $1 \leq i \leq n-1$. In particular, 
      \[
        \varnothing \subset K_{\ADG}(D_{1}) \subset \ldots \subset K_{\ADG}(D_{n-1}) \subset [n]
      \]
      is a strictly increasing nested sequence of sets.\\
      To see $\Phi_{\ADG}(\sigma)=T$, we observe $\sigma^{-1}(1) = K_{\ADG}(D_1)$,
      $\sigma^{-1}(r) = K_{\ADG}(D_{r}) \setminus K_{\ADG}(D_{r-1})$ for $2 \leq i \leq n-1$, and 
      $\sigma^{-1}(n)=[n]\setminus K_{\ADG}(D_{n-1})$. 
      The construction for~$\Phi_{\ADG}(\sigma)$ yields the diagonals $D_{i}$ and the boundary of~$P$,
      in other words, $\Phi(\sigma)=T$.\\
   $(b)\Rightarrow (c)$:
      We trivially have $\{\sigma\} \subseteq (\Phi_{\ADG})^{-1}(T)$, so it remains to show
      $(\Phi_{\ADG})^{-1}(T) \subseteq \{\sigma\}$.

      Assume $\sigma^{\prime} \in S_n$ with $\Phi_{\ADG}(\sigma^{\prime})=T$ and
      $\varnothing \subset K_{\ADG}(D_{1}) \subset \ldots \subset K_{\ADG}(D_{n-1}) \subset [n]$.
      The (unique) singleton set $K_{\ADG}(D_{1}) = \{ r \}$ must consist of a down element: Suppose 
      $r \in K_{\ADG}(D_{1})$ is up, then~$r$ is an endpoint of~$D_1$, otherwise $K_{\ADG}(D_{1}) \neq \{ r\}$.
      But if the other endpoint is~$>r$ (resp.~$<r$) then~$1 \in K_{\ADG}(D_{1})$ (resp. $n \in K_{\ADG}(D_{1})$)
      which also contradicts $K_{\ADG}(D_{1}) = \{ r \}$.
      The diagonal $D_{1}$ must be obtained in the first step of the construction of~$T$ from~$\sigma^{\prime}$, 
      that is, $(\sigma^{\prime})^{-1}(1)=r$, since $r \in K_{\ADG}(D_{j})$ for $j \geq 2$.
      Now suppose, we have finished~$t$ steps in the construction for~$\Phi_{\ADG}(\sigma^{\prime})$ and the 
      diagonals used so far are~$D_{1}, \ldots, D_{t}$. The nestedness of the $K_{\ADG}(D_{i})$ and
      the allowed steps in the construction of~$\Phi_{\ADG}(\sigma^{\prime})$ force 
      $(\sigma^{\prime})^{-1}(t+1) = K_{\ADG}(D_{t+1}) \setminus K_{\ADG}(D_{t})$, i.e. 
      $(\sigma^{\prime})^{-1}(t+1)$ is uniquely determined. Hence $\sigma^{\prime} = \sigma$ by induction.\\
   $(c)\Rightarrow (d)$: 
      Consider $i \in [n]$ with $p_{\ell}^{T}(i)>1$ and $p_{r}^{T}(i)>1$. Then  $2\leq i\leq n-1$.
      Denote by~$u$ the label that realizes~$p_{\ell}^{T}(i)$ and by~$v$ the label that realizes~$p_{r}^{T}(i)$.
      We have $u < i < v$ and $D_u := \{ u,i\}$, $D_v := \{i,v\}$ are diagonals of~$T$. These two diagonals cut~$P$
      into three uniquely determined $\delta_i$-gons~$P_{i}$, $1 \leq i \leq 3$, such that $P_{1}$ is given by~$D_{u}$
      and the path between~$i$ and~$u$ on the boundary of~$P$ that uses only vertices with labels~$\leq i$,~$P_3$ is given
      by~$D_{v}$ and the path between~$i$ and~$v$ on the boundary of~$P$ that uses only vertices with labels~$\geq i$, 
      and~$P_2$ is given by~$D_{u}$, $D_{v}$, and all edges of~$P$ not used by~$P_1$ or~$P_3$.
      The diagonals of~$T$ different from~$D_{u}$ and~$D_{v}$ are a diagonal of precisely one of the polygons~$P_{1}$, 
      $P_{2}$, or $P_{3}$, so we have induced triangulations~$T_i$ of~$P_i$. See Figure~\ref{fig:expl_c_implis_d} for
      an example of this situation.
      \begin{figure}
         \begin{center}
            \begin{minipage}{0.95\linewidth}
               \begin{center}
                  \begin{overpic}
                     [height=6cm]{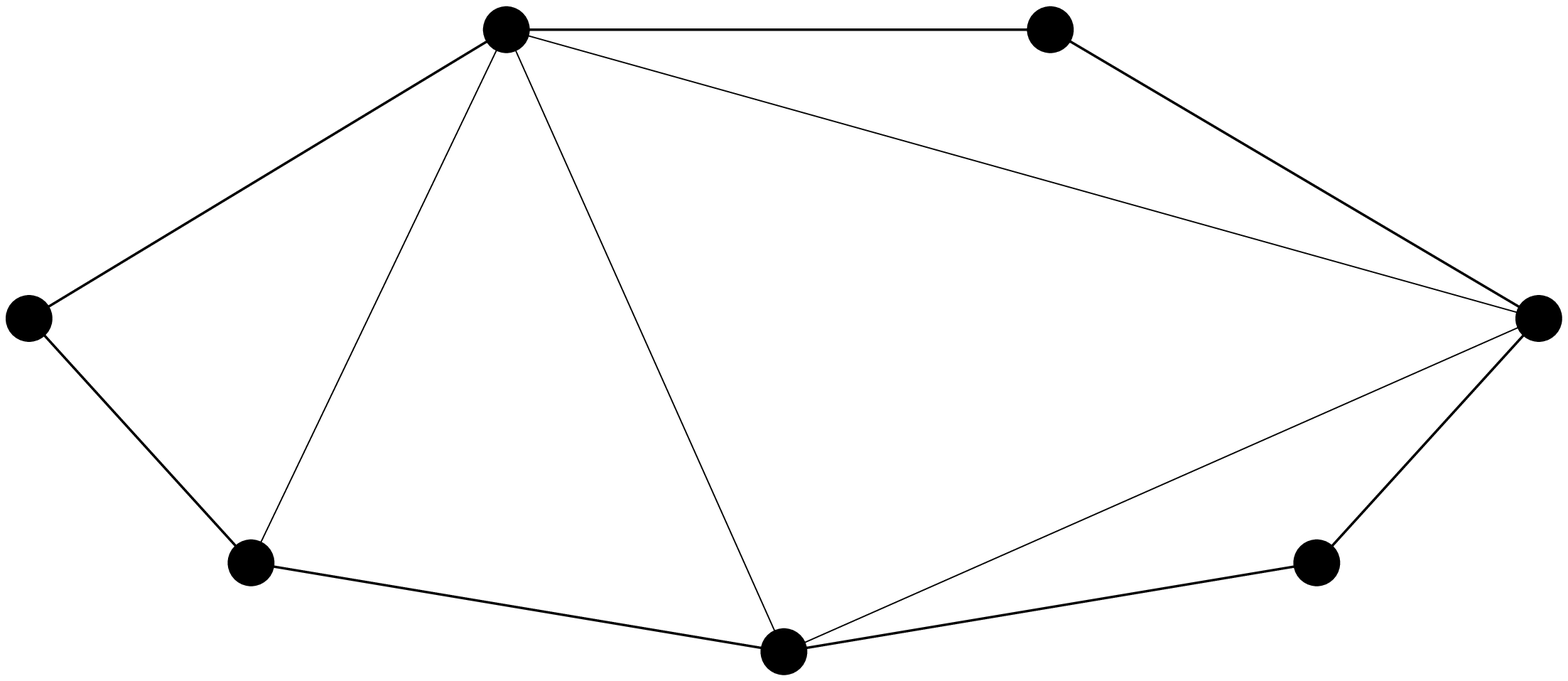}
                     \put(7,28){$0$}
                     \put(22,11){$1$}
                     \put(36,47){$2$}
                     \put(51.5,6){$3$}
                     \put(66.5,47){$4$}
                     \put(81.5,11){$5$}
                     \put(96,28){$6$}
                  \end{overpic}
               \end{center}
               \caption[]{A triangulation~$T$ of a labelled heptagon with $p_{\ell}(3)=3$ and $p_{r}(3)=2$.
                          The polygons~$P_1$,~$P_2$, and $P_3$ have vertex set $\{ 0,1,2,3\}$, $\{ 2,3,4,6\}$,
                          and $\{ 3,5,6\}$. In this case, $\sigma^{\prime} \in S_{5}$ is given by 
                          $(\sigma^{\prime})^{-1} = 21534$, while~$\sigma^{\prime\prime}$ is given by $52134$.}
               \label{fig:expl_c_implis_d}
            \end{minipage}
         \end{center}
      \end{figure}

      Assume that~$i \in \Do_{\ADG}$. There is a (not necessarily unique) permutation~$\sigma^{\prime} \in S_{n}$ with 
      $(\sigma^{\prime})^{-1} = \lambda_1\ldots\lambda_{\delta_1-2}\mu_1\ldots\mu_{\delta_3-2}\nu_1\ldots\nu_{\delta_2-2}$
      such that the $\lambda_i$ generate~$P_1$ and its triangulation~$T_1$, the $\mu_i$ generate~$P_3$ and its triangulation~$T_3$,
      and the $\nu_1$ generate the diagonals of $T_2$ and the path between~$u$ and~$v$ on the boundary of~$P$ that does not contain~$i$.
      But the permutation~$\sigma^{\prime\prime}$ with 
      $(\sigma^{\prime\prime})^{-1} = \mu_1\ldots\mu_{\delta_3-2}\lambda_1\ldots\lambda_{\delta_1-2}\nu_1\ldots\nu_{\delta_2}$
      satisfies $\Phi_{\ADG}(\sigma^{\prime\prime})=T$ and $\sigma^{\prime} \neq \sigma^{\prime\prime}$, that is, 
      $(\Phi_{\ADG})^{-1}(T)$ is not a singleton set.
      
      The case $i \in \Up_{\ADG}$ is handled similarly: Using the same convention
      for~$\lambda_{i}$,~$\mu_{i}$, and~$\nu_{i}$ as in the preceeding case, we find a permutation~$\sigma^{\prime}$ with
      inverse $(\sigma^{\prime})^{-1} = \nu_1\ldots\nu_{\delta_2-2}\mu_1\ldots\mu_{\delta_3-2}\lambda_1\ldots\lambda_{\delta_1-2}$
      which yields~$T$. The permutation~$\sigma^{\prime\prime}$ given by  
      $(\sigma^{\prime\prime})^{-1} = \nu_1\ldots\nu_{\delta_2-2}\lambda_1\ldots\lambda_{\delta_1-2}\mu_1\ldots\mu_{\delta_3-2}$
      is different from~$\sigma^{\prime}$ and satisfies~$\Phi_{\ADG}(\sigma^{\prime\prime})=T$.\\
   $(d)\Rightarrow (a)$: 
      We have to prove that $x_{\sigma^{-1}(i)}=i$, for all $i\in [n]$. Observe
      that $\sigma^{-1}(1)$ must be a down element since $p_{\ell}^{T}(\sigma^{-1}(1))$ or
      $p_{r}^{T}(\sigma^{-1}(1))=1$ and $1,n$ are down elements.
      Therefore $p_{\ell}^{T}(\sigma^{-1}(1))=p_{r}^{T}(\sigma^{-1}(1))=1$ and $x_{\sigma^{-1}(1)}=1$.

      We first prove by induction on $i>1$ that the following assertions are true if $p_{\ell}^{T}(i)=1$ or $p_{r}^{T}(i)=1$.
      \begin{compactenum}[(i.)]
         \item Only one diagonal $D_i=\{u_i,v_i\}$ ($u_i<v_i$) is added at the step $i$ of the
               construction of $\Phi_\ADG(\sigma)$.
         \item The set $\{D_1,\cdots, D_i\}$ defines a triangulation $T_i$ in $\SOT_{i+2}$. The
               set of vertices of this triangulation is then $\{\sigma^{-1}(k),u_k,v_k,\ k\in[i]\}$.
         \item $\{\sigma^{-1}(k),\ k\in[i]\}= K_\ADG (D_i)$.
      \end{compactenum}
      The case~$k=1$ follows from the fact that $\sigma^{-1}(1)$ is a down element. Assume now that these statements
      are true for any $k\in[i-1]$.

      If $\sigma^{-1}(i)\in \Do_\ADG$, then assertion~(i.) holds by definition.
      Statement~(ii.) follows from (c) and by induction (assertion (ii.)), while statement~(iii.) follows
      from statement~(ii.) and by induction (assertion (iii.)).

      Suppose now that $\sigma^{-1}(i)\in \Up_\ADG$. As $p_{\ell}^{T}(\sigma^{-1}(i))$ or
      $p_{r}^{T}(\sigma^{-1}(i))=1$ and by induction (assertion (ii.)), we have to add~$\sigma^{-1}(i)$ to
      the path constructed at the step $i-1$ between the vertices of~$D_{i-1}$, and one of these vertices
      is either a vertex preceding $\sigma^{-1}(i)$ or a vertex following $\sigma^{-1}(i)$ in the labelled
      $(n+2)$-gon. The statements~(i.),~(ii.) and~(iii.) follow now easily from this discussion.

      \medskip
      We now finish the proof: From Lemma~\ref{lem:LeAss1}, Corollary~\ref{lem:LeAss2}, and statement~(iii.) we
      have $\sum_{k\in [i]} x_{\sigma^{-1}(k)} =\tfrac{i(i+1)}{2}$ for all $i\in [n]$. Therefore for $i>1$
      \[
        x_{\sigma^{-1}(i)} = \sum_{k\in [i]} x_{\sigma^{-1}(k)}-\sum_{k\in [i-1]} x_{\sigma^{-1}(k)} = i.
      \]
\end{proof}

\section{The cyclohedron and proofs for Subsection~\ref{sse:IntroCyclo}}\label{sec:RealizCyclo}

In this Section,~$\BDG$ is an orientation of $\BCG_{n}$ or equivalently asymmetric orientation of~$\ACG_{2n-1}$.

\medskip
The convex hull of $\set{M(\sigma)}{\sigma\in W_n}$ is called permutahedron~$\Perm(B_{n})$
of type~$B$. It is well-known that its vertex set is~$\set{M(\sigma)}{\sigma\in W_n}$.
We start this section with an H-representation of~$\Perm(\BCG_n)$. For each $i\in [n]$ we 
consider the hyperplane
\[
  H^{B}_i = \set{x \in \mathbb R^{2n}}{x_{i}+x_{2n+1-i} = 2n+1}.
\]
Such a hyperplane is called {\em type $B$ hyperplane}.  Observe that
$H \cap \bigcap_{i\in[n]} H^{B}_{i}=H\cap \bigcap_{i\in[n-1]} H^{B}_{i}$.
We claim that 
\[
  \Perm(\BCG_{n}) = \bigcap_{i\in [n]} H_{i}^{B} \cap \bigcap_{\varnothing \neq K \subset [n]}\mathscr H_{K}
                  = \Perm(\ACG_{2n-1}) \cap \bigcap_{i\in[n-1]} H^{B}_{i}.
\]
We certainly have $\Perm(B_n) \subseteq \Perm(\ACG_{2n-1}) \cap \bigcap_{i\in[n-1]} H^{B}_{i}$. 
Suppose that~$v \not \in \set{M(\sigma)}{\sigma\in W_n}$ is a vertex
of $\Perm(\ACG_{2n-1}) \cap \bigcap_{i\in[n-1]} H^{B}_{i}$. Then~$v$ must be contained in
the relative interior of an edge of~$\Perm(A_{2n-1})$ which is not entirely contained in $\bigcap_{i\in[n]} H^{B}_{i}$, that is, 
\[
  v = (v_1, \ldots , v_{2n}) =\lambda M(\sigma_1) + (1-\lambda)M(\sigma_2)
\]
for $0 < \lambda <1$ and $\sigma_1,\sigma_2 \in S_{2n}\setminus W_n$ 
with $\sigma_1 = \tau_i\sigma_2$ and  $i\in [2n]\setminus\{n\}$.
If $i<n$, then there is an index $j \in [n]\setminus \{i, i+1\}$ (if $i>n+1$ there is an
index $j \in [n]\setminus \{ 2n+1-i, 2n+2-i\}$) such that 
$\sigma_1(j) + \sigma_1(2n+1-j) \neq 2n+1$. Now $v_j+v_{2n+1-j}\neq 2n+1$
since $\sigma_1(j) = \sigma_2(j)$ and $\sigma_1(2n+1-j) = \sigma_2(2n+1-j)$.
We conclude
\[
  \Perm(\BCG_n) = \operatorname{conv}\set{(\sigma(1),\dots,\sigma(2n))}{\sigma\in W_n\subset S_{2n}}
                = \Perm(\ACG_{2n-1}) \cap \bigcap_{i\in[n-1]} H^{B}_{i}.
\]

\begin{lem}\label{lem:PointsB} 
    Let $T\in \SOT_{2n+2}$. Then $T\in \SOT_{2n+2}^B$ if and only if $M_{\BDG}(T)\in H^{B}_{i}$ 
    for all $i\in [n-1]$.
\end{lem}
\begin{proof} 
   We start with a fundamental observation for any orientation~$\mathscr B$ of~$\BCG_{n}$. As~$\BDG$ is a  
   symmetric orientation of~$\ADG_{2n-1}$, we have $i \in \Do_\BDG\setminus\{1,2n\}$ if and only if 
   $2n+1-i\in\Up_\BDG$ for all $i\in [2n-1]\setminus\{1\}$. A centrally symmetric triangulation~$T$ yields 
   therefore $\omega_i=\omega_{2n+1-i}$, or equivalently,  $x_i+x_{2n+1-i}=2n+1$ for all $i\in\{2,\dots,2n-1\}$. 
   Thus $M_\BDG (T)\in H^B_i$ for any $i\in\{2,\dots,n\}$. From $M_{\BDG}(T)\in H$ follows $x_1+x_{2n}=2n+1$.

   \smallskip
    We now aim for the converse, i.e. consider a triangulation~$T$ of a $(2n+2)$-gon~$P$ with
    $M_{\BDG}(T)\in H^{B}_{i}$ for all $i\in [n-1]$. As $M_\BDG (T)\in H \cap \bigcap_{i\in [n-1]} H_i^B$, we
    conclude that $M_\BDG(T)\in H_n^B$.

    Let us consider a regular realization of the $(2n+2)$-polygon~$P$ labelled according to~$\mathscr B$ 
    and agree on the following terminology: Two labels~$i$ and~$j$ are centrally symmetric if the vertices 
    of~$P$ labelled~$i$ and~$j$ are centrally symmetric. If we consider a triangulation~$T$ of~$P$ then the 
    notion easily extends to edges and triangles. The fundamental observation can now be phrased as label~$i$ 
    is centrally symmetric to label~$2n+1-i$ for $i\in [2n]\setminus\{1,2n\}$. Moreover, label~$1$ is 
    centrally symmetric to label~$2n+1$ and label~$0$ is centrally symmetric to label~$2n$.

    We therefore suppose $\omega_i=\omega_{2n+1-i}$ for $i\in [n]\setminus \{1\}$ and $\omega_1+\omega_{2n}=2n+1$ 
    since $M_\BDG(T)\in H_i^B$ for all $i\in [n]$.

    Choose labels~$a_i$ and~$b_i$ such that $\mu_i(a_i)=p_\ell^T(i)$ and $\mu_i(b_i)=p_r^T(i)$ for all labels $i \in [2n]$. 
    It is easy to see that $\Delta_{i}:=\{ a_{i}, i, b_{i}\}$ is a triangle used by~$T$ to triangulate~$P$.
    As $a_i<i<b_i$, two triangles $\Delta_i$ and $\Delta_j$ coincide if and only if $i=j$. Since~$T$ consists 
    of~$2n$ distinct triangles, the triangles $\Delta_i$ are precisely the triangles used 
    by~$T$. In other words, $\set{\Delta_i}{i\in[2n]}$ determines the triangulation~$T$. 

    We now show by induction on~$k \in [n]$ that~$\Delta_{k}$ is centrally symmetric to~$\Delta_{2n+1-i}$.
    This concludes the proof since~$T$ is a centrally symmetric triangulation of~$P$ if and only if~$\Delta_i$ 
    and $\Delta_{2n+1-i}$ are centrally symmetric for all~$i$.

    If $k=1$ then $a_1=0$ and $b_{2n}=2n+1$. This implies $\mu_1(b_1)=\omega_1=2n+1-\omega_{2n}=2n+1-\mu_{2n}(a_{2n})$. 
    Therefore the edge~$\{a_{2n},2n+1\} \in T$ and $\{1,b_1\} \in T$ are centrally symmetric. Hence the 
    triangles~$\Delta_1$ and~$\Delta_{2n}$ are centrally symmetric.

    Suppose the induction hypothesis is true for $i\in [k]$ where $1<k<n$. If~$a_{k+1}$ does not precede~$k+1$ 
    then $\{a_{k+1},k+1\}$ must be diagonal of~$T$, i.e. an edge of the triangles~$\Delta_{k+1}$ and~$\Delta_\beta$. 
    We conclude from~$a_{k+1}<k+1$ that $k+1=b_\beta$ or $\beta =a_{k+1}$. Both cases imply $\beta\leq k$. In other words,
    there is $\beta\in [k]$ such that $\{a_{k+1},k+1\}$ is an edge of $\Delta_\beta$ or $a_{k+1}$ precedes $k+1$. 

    In the first case $\mu_{k+1}(a_{k+1})=\mu_{2n-k}(b_{2n-k})=:p$ since $\Delta_\beta$ and $\Delta_{2n+1-\beta}$ 
    are centrally symmetric by induction. Hence
    \[
      \mu_{k+1}(b_{k+1})=\tfrac{\omega_{k+1}}{p}=\tfrac{\omega_{2n-k}}{p}=\mu_{2n-k}(a_{2n-k}).
    \]
    Thus $\Delta_{k+1}$ and $\Delta_{2n-k}$ are centrally symmetric. 

    In the second case, the symmetry of $\Delta_k$ and $\Delta_{2n+1-k}$ implies that the label~$b_{2n-k}$ succeeds
    the label~$2n-k$. Again, $\Delta_{k+1}$ and $\Delta_{2n-k}$ are centrally symmetric.
\end{proof}

\medskip
\begin{prop}\label{prop:AssB}
   Let~$\mathscr B$ be an orientation of the Coxeter graph~$B_{n-1}$.
   \begin{compactenum}[1.]
      \item  For~$T\in \SOT^B_{2n+2}$ we have
             \[
               \left\{ M_\BDG(T) \right\}
               = H \cap \bigcap_{D\succ T} H_{K_{\BDG}(D)} \cap \bigcap_{i\in[n-1]} H^{B}_{i}.
             \]
      \item The intersection of the hyperplane~$H$, the type~$B$ hyperplanes~$H^{B}_{i}$, $i\in[n-1]$,
            and the $\BDG$-admissible half spaces~$\mathscr H_{K}$ is an \textnormal{H}-representation
            of the cyclohedron $\Ass(\BCG_{n})$. In particular, the
            permutahedron~$\Perm(\BCG_n)$ is contained in the cyclohedron $\Ass(\BCG_n)$ which
            is contained in the associahedron~$\Ass(\ACG_{n-1})$.
   \end{compactenum}
\end{prop}

\begin{proof}
\begin{compactenum}[1.]
  \item follows from Lemma~\ref{lem:PointsB} and Corollary~\ref{cor:AssAPoints}.

  \item We first observe that the intersection of all admissible half spaces and of all  type~$B$ hyperplanes
  defines a bounded
   set in~$\R^{2n}$. This follows immediately from Theorem~\ref{cor:AssABound}.
   The intersection of all those admissible hyperplanes with~$H$ yields therefore a bounded $n$-dimensional
   convex polytope.

   The first part, Corollary~\ref{cor:AssAPoints}, Theorem~\ref{cor:AssABound}, and Lemma~\ref{lem:PointsB} show
   that the set of vertices of this polytope is $\{M_\BDG(T),\, T\in \SOT^B_{2n+2}\}$ and that this convex 
   polytope is simple: each vertex is contained in precisely $2n-1 - (n-1) = n$ facet defining hyperplanes.

   A centrally bistellar flip in a centrally symmetric triangulation is a succession of at most
   two bistellar flips: flip a diagonal together with its centrally symmetric.
   By (1.), the $1$-skeleton of this polytope is the flip graph of the centrally symmetric triangulations
   of an $(2n+2)$-gon: Two vertices `differ in precisely one centrally diagonal flip' if and only if the
   vertices are connected by an edge. Therefore it is the
   $1$-skeleton of the cyclohedron (\cite[Theorem 1]{simion}).
\end{compactenum}
\end{proof}

\subsection{Proof of Theorem~\ref{thm:AssB}, Proposition~\ref{prop:AssBfromPermB}, and Proposition~\ref{prop:AssBandPermB}} 
$ $\\
These statements follow immediatley from Propositions~\ref{prop:AssAfromPermA},~\ref{prop:AssAandPermA}, and~\ref{prop:AssB}.

\section{Remarks and Questions}\label{sec:Remarks}

\subsection{On normal fans of these realizations}\label{sse:cambrian}

Recall the following well-known facts about
Coxeter groups. Let $W$ be a finite Coxeter group acting on a
vector space $V$ as a reflection group. The Coxeter fan of $W$
(relative to $V$) is the fan created by the Coxeter (hyperplane)
arrangement of $W$ in $V$. Choose a generic point in a maximal cone of
the Coxeter fan, then the convex hull of the orbit of this point
under the action of $W$ yields a permutahedron whose normal fan is
the Coxeter fan. An example are the permutahedra~$\Perm(\ACG_{n-1})$
and~$\Perm(\BCG_n)$ which are the convex hull of the $S_n$-orbit 
of~$(1,2,\ldots,n)$ and the $W_n$-orbit of~$(1,2,\ldots,2n)$.
Denote by~$\mathcal N(\ACG_{n-1})$ the  normal fan of~$\Perm(\ACG_{n-1})$. 
This fan is a Coxeter fan of type $\ACG_{n-1}$ in $H$.  

For each orientation~$\ADG$ of~$\ACG_{n-1}$, the {\em cambrian fan}~$\mathcal N(\ADG)$ 
associated to the orientation~$\ADG$ is the fan obtained by gluing all maximal cones 
in~$\mathcal N(\ACG_{n-1})$ that correspond to permutations~$\sigma \in \Phi_\ADG^{-1}(T)$ 
for any given~$T \in \SOT_{n+2}$, \cite{reading}. Reading proved that this fan is always
simplicial and the normal fan of a realization of the associahedron that comes from Loday's 
or a bipartite orientation. He conjectured that this is true for any orientation 
of the Coxeter graphs of type~$A$ and~$B$.

Moreover, Reading explicitly described the rays of~$\mathcal N(\ADG)$ in~\cite[Sec.~9]{reading}.
This description translates directly to a description of the admissible half spaces
for a given orientation~$\ADG$. It is easy to see that this bijection extends to
an isomorphism of the face lattice of~$\mathcal N(\ADG)$ and the normal fan of the 
associahedron obtained from~$\ADG$. Hence we obtain the following proposition.

\begin{prop}\label{cor:AssAHrep}
   Fix an orientation~$\ADG$ on~$\ACG_{n-1}$. The normal fan
   of the realization of $\Ass(\ACG_{n-1})$ associated to $\ADG$ is precisely the
   cambrian fan $\mathcal N (\ADG)$.
\end{prop}

The normal fan of the type $B$ permutahedron
$\Perm(B_n)$ lives in  $\mathbb R^{2n}\cap  H\bigcap_{i\in[n-1]}
H^B_i$ and is precisely
$$
\mathcal N(\BCG_n)=\mathcal N(\ACG_{2n-1})\cap \bigcap_{i\in[n-1]}
H^B_i.
$$
Let $\BDG$ be a symmetric orientation on $\ACG_{2n-1}$ or
equivalently an orientation of $\BCG_n$. From Reading's work (last
sentence of \cite{reading}): The Cambrian fan $\mathcal N_B(\BDG)$
 of type $B_n$ is given from the
corresponding cambrian fan $\mathcal N(\BDG)$ of type $A_{2n-1}$
by the formula:
\[
  \mathcal N_B(\BDG)=\mathcal N(\BDG)\cap \bigcap_{i\in[n-1]} H^B_i.
\]
Hence we have the following corollary..

\begin{cor}
  Fix an orientation~$\BDG$ on~$\BCG_{n}$. The normal fan
  of the realization of $\Ass(\BCG_{n})$ associated to $\BDG$ is  the
  cambrian fan $\mathcal N_B (\BDG)$.
\end{cor}

\begin{rem}
\textnormal{
In \cite[\S9]{reading}, N.~Reading proved in type $A_{n}$ and
$B_{n}$ that the cambrian fan corresponding to a bipartite
orientation (i.e. $i$ is a down element if and only if $i+1$ is an
up element) is linearly isomorph to a cluster fan. The
realization of the permutahedron of type $A_n$ or $B_n$ used in
this article fix a {\em geometric representation} of the corresponding
Coxeter group. Let $\Phi$ be a crystallographic  root system and
$\Phi^+$ be its set of positive roots: the cluster fan associated
to $\Phi$ is then the fan spanned by the almost positive roots of
$\Phi$ \cite{fomin_zelevinsky}. Hence, for a bipartite orientation, we
have a realization of the associahedron (or of the cyclohedron)
whose normal fan is linearly isomorph to a cluster fan, as in
\cite{chapoton_fomin_zelevinsky} (see also \cite[Theorem 5.11]{fomin-reading}).
}
\end{rem}

\subsection{On isometry classes of these realizations}\label{sse:Isometry}

We are starting here a study of the (affine) isometry classes of our
realizations of $\Ass(\ACG_{n-1})$. Some experiments with
{\tt GAP}~\cite{gap} and {\tt polymake}~\cite{polymake} show that these realizations
are not all isometric. Indeed, if two realizations are isometric,
then they necessarily have the same number of common vertices with
$\Perm(\ACG_{n-1})$ but we shall see that this condition is not sufficient.
It would be interesting to classify the isometry
classes of these realizations in terms of `equivalence classes'
on orientations of the Coxeter graph~$\ACG_{n-1}$.

The simplest definition of such equivalence classes yields
isometric realizations. Two orientations~$\ADG$
and~$\ADG^{\prime}$ of~$\ACG_{n-1}$ are {\em equivalent} if~$\ADG$
is obtained from~$\ADG^{\prime}$ by reversing the orientations of
all edges. This implies
$\Up_{\ADG}=\Do_{\ADG^{\prime}}\setminus\{1,n\}$ and each
equivalence class consists of two orientations. The following
result can be easily deduced from definitions:

\begin{prop}\label{prop:isometry}
   Let $\ADG$ and $\ADG^{\prime}$ be two orientations of $\ACG_{n-1}$. If $\ADG$ and $\ADG^{\prime}$
   are equivalent,  then the isometric transformation $(x_1,\dots,x_n) \to (n+1-x_1,\dots,n+1-x_n)$
   on $\mathbb R^n$ maps the realization of $\Ass(\ACG_{n-1})$ associated to $\ADG$ on the realization
   of $\Ass(\ACG_{n-1})$ associated to $\ADG^{\prime}$.
\end{prop}

Each orientation is completely determined by its set of up indices. The following table gives the
number~$n_\ADG$ of common vertices of~$\Perm(\ACG_{n-1})$ and~$\Ass(\ACG_{n-1})$ for each orientation~$\ADG$
of~$\ACG_{n-1}$ for~$n\leq 5$ as well as the number $I_{\ADG}$ of integer points contained in the
associahedron. The number~$n_{\ADG}$ can be either computed by {\tt GAP}, with an algorithm based on
the cambrian congruences, and the equivalence between $(a)$ and $(c)$ in Proposition~\ref{prop:AssAandPermA}
or by counting the vertices of the associahedron with coordinates a permutation of~$(1,2,\ldots,n)$. The
coordinates can be obtained for example by using {\tt polymake}, the numbers~$I_{\ADG}$ can be computed 
with the help of {\tt LattE},~\cite{latte}. Input data for all examples is available at~\cite{lange}. 
{\scriptsize
\[
\begin{array}{l|cc|cccc|cccccccc}
         & n=3        &     &n=4        &       &     &     &n=5        &         &     &       &     &       &       &     \\\hline
\Up_\ADG & \varnothing&\{2\}&\varnothing&\{2,3\}&\{2\}&\{3\}&\varnothing&\{2,3,4\}&\{2\}&\{3,4\}&\{4\}&\{2,3\}&\{2,4\}&\{3\}\\
n_\ADG   & 4          &4    &8          &8      &9    &9    &16         &16       &19   &19     &19   &19     &20     &20   \\
I_{\ADG} & 8          &8    &55         &55     &60   &60   &567        &567      &672  &672    &672  &672    &742    &742
\end{array}
\]
}
For $n=5$, $\{2\}$ and $\{3,4\}$ form an equivalence class as well as $\{4\}$ and $\{2,3\}$. All these up sets yield
$n_{\ADG}=19$, and the number of integer point they contain is~$672$. 

We now consider the transitive closure of the following modification of the equivalence of two orientations introduced above.
This modified notion of equivalence yields the equivalence classes $\varnothing$ and $\{2,3,4\}$; $\{2\}$, $\{3,4\}$,
$\{4\}$, and $\{2,3\}$; and $\{ 2,4\}$ and $\{3\}$ in case of $n=5$. 
Two orientations~$\ADG$ and~$\ADG^{\prime}$ of~$\ADG_{n-1}$ are {\em equivalent} if~$\ADG$ is obtained 
from~$\ADG^{\prime}$ by reversing the orientations of all edges or if the oriented graph~$\ADG$
is obtained from ~$\ADG^{\prime}$ by a rotation of 180 degrees. The transitive closure of this modified 
notion of equivalence suggests isometry classes for $n=6$ that can be detected by~$I_{\ADG}$
but not by~$n_{\ADG}$:
\[
\begin{array}{l|cccccccc}
   n=6        &             &           &           &         &           &           &         &        \\ \hline
   \Up_{\ADG} & \varnothing & \{ 2 \}   & \{5\}     & \{2,3\} & \{3\}     & \{4\}     & \{2,5\} & \{2,4\}\\
              & \{2,3,4,5\} & \{3,4,5\} & \{2,3,4\} & \{4,5\} & \{2,4,5\} & \{2,3,5\} & \{3,4\} & \{3,5\}\\
   n_{\ADG}   &     32      &    39     &    39     &    42   &    42     &    42     &    44   &    45  \\
   I_{\ADG}   &    7958     &   10116   &   10116   &   11155 &   12294   &   12294   &   12310 &  13795 
\end{array}
\]
We believe that these equivalence classes can be characterized by the number of integer lattice points contained
by the corresponding realizations, and that such an equivalence class consists precisely of isometric realizations.
Moreover, viewing these polytopes as generalized permutahedra as defined by A.~Postnikov~\cite{postnikov}, it should 
be possible to describe these realizations explitly as a Minkowski sum (with possibly negative coefficients). Once a 
such a Minkowski sum decomposition is determined for a given oriented Coxeter graph, explicit formulae for the volume 
and number of integer points are explicit.

\subsection{On barycenters}\label{sse:Barycenters}

In his article, J.-L.~Loday mentions an observation made by
F.~Chapoton that the vertices of the permutahedron and the
associahedron of his realization have the same barycenter:
$G=(\frac{n+1}{2},\ldots,\frac{n+1}{2})$. We observed that for
$n \leq 10$ and any orientation $\ADG$ of $\ACG_{n-1}$, the 
barycenter of the vertices of the realization of the associahedron 
associated to $\ADG$ is $G$. This seems to be true for the 
cyclohedron, too, as we computed all barycenters of cyclohedra up to 
dimension~$5$. This leads us to the following question: 
Let $\ADG$ be an orientation of $\ACG_{n-1}$ and $\BDG$ an orientation of~$\BCG_{n}$, 
is $G$ the barycenter of $\operatorname{conv}\{M_\ADG(T), \ T\in \SOT_{n+2}\}$
and $\operatorname{conv}\{M_\BDG(T), \ T\in \SOT^{B}_{2n+2}\}$?

\bigskip

\subsection*{Acknowledgements}

The authors are grateful to the organizers Anders Bj\"orner and
Richard Stanley of the algebraic combinatorics session at the
Institut Mittag-Leffler in Djursholm, Sweden, where the main part
of this work was done. They also thank  Jean-Louis Loday for some useful
comments on a preliminary version of this work.

C.H.: I wish to thank Jean-Louis Loday for very instructive
conversations on the topic of this article, which took place when
I was a member of the Institut de Recherche Math\'ematique Avanc\'ee in
Strasbourg, France.

C.L.: I thank the organizers Ezra Miller, Vic Reiner, and Bernd Sturmfels of the PCMI summer session 2004 in Park City,
Utah, for making my participation possible. My interest in generalized associahedra was initiated by Sergey Fomin's
lectures {\em Root systems and generalized associahedra} during this session. Moreover, many thanks to the polymake
team (Ewgenij Gawrilow, Michael Joswig, Thilo Schr\"oder, and Niko Witte) for their invaluable help, to Konrad Polthier
(the visualization would be much harder without javaview~\cite{javaview}), and to Jes\'us de Loera for making Latte available.


\end{document}